\documentclass[10pt]{amsart}

\theoremstyle{plain}
\newtheorem{thm}{Theorem}[section]

\newtheorem{lem}[thm]{Lemma}
\newtheorem{prop}[thm]{Proposition}
\newtheorem{cor}[thm]{Corollary}

\theoremstyle{definition}

\newtheorem{df}[thm]{Definition}
\newtheorem{dfs}[thm]{Definitions}
\newtheorem{ex}[thm]{Example}
\newtheorem{ex-notn}[thm]{Example/Notation}

\newtheorem{conj}[thm]{Conjecture}
\newtheorem{prob}[thm]{Problem}
\newtheorem{rem}[thm]{Remark}

\newtheorem{notn}[thm]{Notation}

\def\ann{\operatorname{ann}}

\def\ddeg{\operatorname{-deg}}
\def\der{\operatorname{der}}
\def\Div{\operatorname{div}}
\def\DR{\operatorname{DR}}

\def\gr{\operatorname{gr}}

\def\im{\operatorname{im}}

\def\var{\operatorname{Var}}

\def\del{\partial}

\def\calA{{\mathcal A}}
\def\calM{{\mathcal M}}
\def\calN{{\mathcal N}}
\def\calP{{\mathcal P}}
\def\CC{{\mathbb C}}

\def\calD{{\mathcal D}}
\def\euler{{\mathcal E}}
\def\fraka{{\mathfrak a}}
\def\frakm{{\mathfrak m}}
\def\frakn{{\mathfrak n}}
\def\frakA{{\mathfrak A}}

\def\NN{{\mathbb N}}

\def\QQ{{\mathbb Q}}    

\def\ZZ{{\mathbb Z}}

\def\action{\bullet}
\def\divides{\operatorname{|}}
\def\dividesnot{{\not |\,\,}}
\def\bar#1{\overline{#1}}

\def\into{\hookrightarrow}
\def\Mtwo{{\em Macaulay} 2\expandafter}
\def\onto{\to\hskip-1.7ex\to}

\def\ideal#1{\langle #1 \rangle}

\def\mylabel#1{\label{#1}}
\def\ignore#1{}

\numberwithin{equation}{section}
\begin{document}

\title[Bernstein-Sato polynomial and Milnor fiber cohomology]
       {Bernstein-Sato polynomial versus cohomology of the Milnor
       fiber for generic hyperplane arrangements}
\thanks{The author wishes to thank the National Science Foundation,
the Deutsche Forschungsgemeinschaft, and the
    A.v.-Humboldt-Stiftung for their support.}
\author[Uli Walther]{Uli Walther}
\address{Department of Mathematics,
Purdue University,
150 North University Street,
West Lafayette, IN  47907-2067}
\email{walther@math.purdue.edu}
\begin{abstract}
Let $Q\in\CC[x_1,\ldots,x_n]$ be a homogeneous polynomial
of degree $k>0$. We 
establish a connection between the Bernstein-Sato polynomial $b_Q(s)$ 
and the degrees of the generators for the top cohomology of the
associated Milnor fiber. 
In particular, the integer $u_Q=\max\{i\in\ZZ:
b_Q(-(i+n)/k)=0\}$ bounds the top degree (as differential form) of the
elements in $H^{n-1}_{\DR}(Q^{-1}(1),\CC)$. 
The link is provided by the relative de Rham
complex and $\calD$-module algorithms for computing integration
functors. 

As an application we determine the Bernstein-Sato polynomial $b_Q(s)$
of a generic  
central arrangement $Q=\prod_{i=1}^kH_i$ of hyperplanes. We obtain
in turn information about the cohomology of the Milnor fiber of such
arrangements related to results of Orlik and Randell who investigated
the monodromy.

We also introduce certain
subschemes of the arrangement determined by the roots of
$b_Q(s)$. They appear to correspond to iterated singular loci.
\end{abstract}
\maketitle

\section{Introduction}
Let $f$ be a non-constant polynomial in $n$ variables. 
In the 1960s, M.~Sato introduced $a$-, $b$- and
$c$-functions associated to a prehomogeneous vector space
\cite{Satoetal,Sato}.  The existence of $b$-functions associated to
all polynomials 
and germs of holomorphic functions was later established in
\cite{Bernstein,Bjork2}. 

The simplest interesting case of a $b$-function is the
case of the quadratic form $f(x_1,\ldots,x_n)=\sum_{i=1}^nx_i^2$. Let
$s$ be a new variable and denote by $f^s$ the germ of the complex
power of $f(x)$. One
then has an identity 
\[
\left(\sum_{i=1}^n\frac{\del^2}{\del {x_i}^2}\right)
\action f^{s+1}=4(s+1)(s+n/2)f^s.
\]
The $b$-function to $f(x)$ is here
$b_f(s)=(s+1)(s+n/2)$. One may for general $f$ 
use an equality of the type
\begin{eqnarray}
\label{eqn-bpoly-def}
P(s)\action f^{s+1}&=&b(s)f^s
\end{eqnarray}
to analytically
continue $f^s$, and it was this application that initially caused
I.N.~Bernstein to consider $b_f(s)$. 
Today, the $b$-function
of a polynomial is usually referred to as ``Bernstein-Sato
polynomial'' and denoted $b_f(s)$. 

The Bernstein-Sato polynomial
is always a multiple of $(s+1)$, and equality holds if $f$ is smooth. 
The roots of $b_f(s)$ are always negative and rational
\cite{Kashiwara}. 
It
has been pointed out first in \cite{Malgrange2,Malgrange} that there
is an intimate connection between 
the singularity structure of $f^{-1}(0)$
and its Bernstein-Sato polynomial.
The roots of $b_f(s)$ relate to a variety of
algebro-geometric data like the structure of the embedded
resolution of the pair $(\CC^n,\var(f))$, Newton polyhedra, 
Zeta functions, asymptotic expansions of integrals, Picard-Lefschetz
monodromy, polar invariants and multiplier ideals: see, for example, 
\cite{Budur-Saito,Hamm2,Kato,Lichtin,Loeser,Varchenko2}. T.~Yano
systematically worked
out a number of examples \cite{Yano} and some interesting computations
are given in \cite{BGMM}. A satisfactory
interpretation of all roots of $b_f(s)$ for general $f$ remains,
however, outstanding. Indeed, until \cite{Oa2} there was not even 
an algorithm for the computation of $b_f(s)$ for an arbitrary
polynomial $f$. 

\bigskip

In this note we investigate the Bernstein-Sato
polynomial when $f$ defines a generic central hyperplane
arrangement. By that we mean a reduced collection of $k$ hyperplanes 
such that each subset of $\min\{k,n\}$ of the hyperplanes
cuts out the origin.
The paper is organized as follows. In this section we introduce the
relevant notation. In the next section 
we find an upper bound for the Bernstein-Sato polynomial of a central generic
arrangement. We shall compute a polynomial $b(s)$ that satisfies an
identity of the type (\ref{eqn-bpoly-def})  
 using strongly that the arrangement is central and generic.
In Section 3 we use some counting and Gr\"obner type 
arguments to obtain
information about generators for the top cohomology of the Milnor
fiber of such arrangements. We prove parts of a conjecture of Orlik and
Randell on the cohomology of the Milnor fiber of a generic
central arrangement. 
In particular, we
determine in exactly 
which degrees the top cohomology lives, and we present a
conjectured set of generators. 

Malgrange \cite{Malgrange3} 
demonstrated that the Bernstein-Sato polynomial is the
minimal polynomial of a certain operator on the sheaf of vanishing
cycles. This says in essence that monodromy eigenvalues are exponentials
of roots of $b_f(s)$. In the fourth section we prove roughly that
for homogeneous $f$ the degrees of the top Milnor fiber
cohomology are roots of $b_f(s)$. This can in some sense 
be seen as a logarithmic
lift of Malgrange's results. For generic central
arrangements this links our results from Sections 2 and 3 and allows
the determination of all roots of $b_f(s)$ and (almost) all
multiplicities. 
We close Section 4 with an example of
 a non-generic arrangement, and finish in Section 5 with some
statements and conjectures 
regarding  the structure of the $D_n$-modules $R_n[f^{-1}]$ and
$D_n[s]\action f^s$.

\begin{notn} 
Throughout, we will work over the field of complex numbers
$\CC$. We should point out that this is mostly for keeping things
simple as the Bernstein-Sato polynomial is invariant
under field extensions.

In this note, for elements $\{f_1,\ldots,f_k\}$ of any ring $A$, 
$\ideal{f_1,\ldots,f_k}$ denotes the left ideal generated by
$\{f_1,\ldots,f_k\}$. If we mean a right ideal, we specify  it by
writing $\ideal{f_1,\ldots,f_k}A$.

By $R_n$ we denote the ring of polynomials $\CC[x_1,\ldots,x_n]$
in $n$ variables over $\CC$, and by $D_n$ 
we mean the ring of $\CC$-linear
differential operators on $R_n$, the $n$-th Weyl algebra. The ring
$D_n$ is generated by the partial derivative operators
$\del_i=\frac{\del}{\del x_i}$ and
the multiplication operators $x_i$. One may consider $R_n$ as a subring
of $D_n$ as well as a quotient of $D_n$ (by the left ideal
$\ideal{\del_1,\ldots,\del_n}$). We denote by $\action$ the natural
action of $D_n$ on $R_n$ via this quotient map, as well as induced
actions of $D_n$ on localizations of $R_n$.

We will have occasion to consider $D_t$, $D_x$ and $D_{x,t}$ in some
instances, where $D_t$ is the Weyl algebra in the variable $t$, $D_x$
the one in  
$x_1,\ldots,x_n$ and $D_{x,t}$ is the Weyl algebra in $x_1,\ldots,x_n$
and $t$. 

The module of global algebraic differential $n$-forms on $\CC^n$ 
is denoted $\Omega$; it
may be pictured as the quotient
$D_n/\ideal{\del_1,\ldots,\del_n}D_n$. The left $D_n$-Koszul complex
on $D_n$ induced by 
the commuting vector fields 
$\del_1,\ldots,\del_n$ is denoted $\Omega^\bullet$; it is a resolution
for $\Omega$ as right $D_n$-module.

We shall use multi-index notation in $R_n$: writing 
$x^\alpha$ implies that $\alpha=(\alpha_1,\ldots,\alpha_n)$ and stands
for $x^\alpha=x_1^{\alpha_1}\cdots x_n^{\alpha_n}$. The same applies to
elements of $D_n$, both for the polynomial and the differential
components. If $\alpha$ is a multi-index, $|\alpha|$ denotes the sum of
its components; if $I$ is a set, then $|I|$ is its
cardinality. Finally, if $k,r\in\NN$ then $k\divides r$ signifies that $k$
divides $r$ while $k\dividesnot r$ indicates that this is not the case.

\end{notn}


\subsection{Bernstein-Sato polynomials}\hfill

\begin{df}
For $f\in R_n$ we define $J(f^s)\subseteq D_n[s]$ to be the
annihilator of $f^s$ via formal differentiation, 
this is a left ideal. We set
\[
\calM=D_n[s]/(J(f^s)+\ideal{f})=D_n\action f^s/D_n\action f^{s+1}.
\] 
By definition, the
{\em Bernstein-Sato polynomial} $b_f(s)$ of $f$ is the minimal
polynomial of $s$ on $\calM$.
So $b_f(s)$ is the monic
polynomial of smallest degree satisfying a functional equation of the
type (\ref{eqn-bpoly-def})
with $P(s)\in D_n[s]$.
\end{df}

Let 
$\tilde\calM=D_n[s]/(J(f^s)+\ideal{f}+D_n[s]\cdot \frakA)$ 
where $\frakA\subseteq
R_n$ is the Jacobian ideal of $f$,
$\frakA=\sum_{i=1}^nR_n\del_i\action(f)$. 
Then $\tilde \calM$
is isomorphic to 
$(s+1)\calM$ and since $(s+1)$ divides $b_f(s)$ then 
the minimal polynomial of $s$ on $\tilde\calM$
is $\tilde b_f(s)=b_f(s)/(s+1)$. 

\bigskip

Consider the module $D_n\action f^a$ for $a\in \CC$ and write $J(f^a)$
for the kernel of the map $D_n\to D_n\action f^a$ induced by $P\mapsto
P\action f^a$. There is a natural map $D_n\action f^{a+1}\into
D_n\action f^{a}$ induced by $P\action f^{a+1}\mapsto Pf\action f^a$. 
Some roots of the
Bernstein-Sato polynomial detect the failure of this map to be an
isomorphism:
\begin{lem}
\mylabel{lem-broot}
Let $a\in\QQ $ be  such that $b_f(a)=0$ but $b_f(a-n)\not =0$
for all positive natural numbers $n$. Then $D_n\action f^a\not =
D_n\action f^{a+1}$. 
\proof
Suppose that $a$ is as the hypotheses stipulate, and in addition
assume that $D_n\action f^a =
D_n\action f^{a+1}$. We will exhibit a contradiction.

Since $D_n\action f^{a+1}\to 
D_n\action f^{a}$ is an epimorphism, $D_n=\ideal{f}+J(f^a)$. By the
choice of $a$ and Proposition 6.2 in \cite{Kashiwara},
$J(f^a)=D_n\cap(J(f^s)+D_n[s]\cdot(s-a))$. Hence
$D_n[s]=J(f^s)+\ideal{f}+\ideal{s-a}$. Multiplying by $b_f(s)/(s-a)$
we get $\ideal{b_f(s)/(s-a)}\subseteq
J(f^s)+\ideal{f}+\ideal{b_f(s)}$. Since $b_f(s)\in J(f^s)+\ideal{f}$, 
\[
\frac{b_f(s)}{s-a}\in J(f^s)+\ideal{f}.
\]
That, however, contradicts the definition of $b_f(s)$ as the minimal
polynomial in $s$ contained in the sum on the right.
\qed
\end{lem}

\subsection{Isolated Singularities}\hfill

Suppose that $f$ has an isolated singularity and assume for simplicity
that the singularity is at the origin. We give a short overview of
what is known about the Bernstein-Sato polynomial in this case,
following \cite{Kashi-book,Malgrange,Yano}.

The module $\tilde \calM$ is supported only at the origin, so 
by \cite{Kashiwara} the minimal
polynomial of $s$ on $\Omega\otimes_{D_n}\tilde\calM$ is
$\tilde b_f(s)$. 
If now $f$ is homogeneous of degree $k$, 
$kf=\sum_{i=1}^n x_i\del_i\action(f)$. 
Then $J(f^s)$ contains $\sum_{i=1}^nx_i\del_i-ks$.
The action of $s$ on
a homogeneous $g\in \Omega\otimes_{D_n}\tilde\calM\cong R_n/\frakA$ is
easily seen to be multiplication by 
$(-n-\deg(g))/k$. Thus, the Bernstein-Sato
polynomial of a homogeneous isolated singularity 
 encodes exactly the degrees of non-vanishing
elements in $R_n/\frakA$. 

Consider now the relative de Rham complex $\Omega_f^\bullet$ 
associated to the map
$f:\CC^n\to \CC$. We shall denote the coordinate on $\CC$ by $t$. The
complex $\Omega^\bullet_f$ is the Koszul complex induced by left
multiplication by $\del_1,\ldots,\del_n$ on the $D_{x,t}$-module
$\calN=D_{x,t}/J_{n+1}(f)$ where $J_{n+1}(f)$ is the left ideal of $D_{x,t}$
generated by $t-f$ and the expressions $\del_i+\del_i\action(f)
\del_t$ for $1\le
i\le n$. 
The complex 
$\Omega_f^\bullet=\Omega^\bullet\otimes_{D_n}\calN$ is a representative
of the application of the de Rham functor $\int_f$ 
associated to the map $f$ to
the structure sheaf on $\CC^n$, \cite{Deligne}.
Its last nonzero cohomology module appears in
degree $n$,
$H^{n}(\Omega_f^\bullet)=\calN/\{\del_1,\ldots,\del_n\}\cdot\calN$.
This module is in a natural way a left $D_t$-module. For any $\alpha\in
\CC$, 
the cohomology of the 
derived tensor product of 
$\Omega^\bullet_f$ 
with
$D_t/\ideal{t-\alpha}D_t$ 
is the de Rham
cohomology of the fiber at $\alpha$. The identification of
$\calN/\{\del_1,\ldots,\del_n,t-\alpha\}D_{x,t}$ with
$H^{n-1}_{\DR}(\var(f-\alpha))$ is explained in and before Lemma
\ref{lem-omega}.

So one has an isomorphism  
\[
R_n/\frakA\cong
(D_t/\ideal{t-\alpha}D_t)\otimes_{D_t}
H^n(\Omega^\bullet\otimes_{D_n}\calN)\cong  
H^{n-1}_{\DR}(f^{-1}(\alpha),\CC)
\]
and the roots of $b_f(s)$ in fact represent the degrees
of the cohomology classes of the Milnor fiber of $f$. 

For general $f$, the Bernstein-Sato polynomial is more
complex,
see Example \ref{ex-non-generic} and the following remarks.

\section{An upper bound for the Bernstein polynomial}

Our goal is Theorem \ref{thm-upper-bound}.
We shall mimic some of the
mechanism that makes the isolated singularity case so easy. It is
clear that a literal translation is not possible, because $R_n/\frakA$
has in general elements in infinitely many different degrees. 
However, we now introduce
 certain ideals in
$R_n$ that are intimately related to the Bernstein-Sato polynomial.

\begin{df}
Let $q(s)\in \CC[s]$. For a fixed $f\in R_n$ 
we define the ideal $\fraka_{q(s)} \subseteq
R_n$ as the set of elements $g\in R_n$ 
\[
\left[g\in \fraka_{q(s)}\right]
\Longleftrightarrow
\left[\exists 
P(s)\in D_n[s]: P(s)\action f^{s+1}=q(s)gf^s\right].
\]
We remark that $\fraka_{q(s)}\subseteq
\fraka_{q(s)q'(s)}$, and if $q'(s)g\in \fraka_{q(s)}\cdot R_n[s]$ then
$g\in\fraka_{q(s)q'(s)}$. The Jacobian ideal $\frakA$ is
contained in $\fraka_{(s+1)}$, and $f\in\fraka_{(1)}$.

The Bernstein-Sato polynomial of $f$ is evidently the polynomial 
$b_f(s)$
of smallest degree
 such that $1\in\fraka_{b_f(s)}$. 
\end{df}
Before we come to the computation of an estimate for $b_f(s)$ for
generic arrangements we first consider general homogeneous polynomials
and then arrangements in the plane.

\subsection{The homogeneous case}\hfill

Assume now that $Q\in R_n$ is
homogeneous.\footnote{Throughout we use $Q$ for an instance of a
  homogeneous polynomial while $f$ is used if no homogeneity
  assumptions are in force.} We shall denote by $\frakm$ the
homogeneous maximal ideal of $R_n$.
If $g\in\fraka_{q(s)}$ then by definition $gQ^s\in\calM$ is annihilated by
$q(s)$. Since $b_Q(s)$ annihilates all of $\calM$, finding
$g\in\fraka_{q(s)}$ is equivalent to finding eigenvectors of
$s$ on $\calM$ to eigenvalues that are zeros of $q(s)$.
In the isolated singularity case one only has to study the
residues of $\fraka_{q(s)}$ in $R_n/\frakA$, and this goes as follows.
Let $\delta_Q
=\min_{k\in\NN}\{\frakm^{k+1}\subseteq  \frakA\}$. 
Then
 the homogeneous polynomial 
$g$ with $0\not = \bar g\in R_n/\frakA$ is in $\fraka_{q(s)}$ if and only if
$(s+1)\prod_{i=\deg(g)}^{\delta_Q}\left(s+\frac{i+n}{\deg(Q)}\right)$ 
divides $q(s)$; this is proved in  \cite{Yano} based on results of
Kashiwara.  

For non-isolated homogeneous singularities $Q$ we have a weak version of this:

\begin{lem}
\mylabel{lem-euler}
If $R_n[s]\cdot 
\fraka_{q(s)}$ contains $\frakm^rg$ where $g=g(s)\in R_n[s]$ is 
homogeneous in $x_1,\ldots,x_n$  
then $g\in
R_n[s]\cdot \fraka_{q'(s)}$ where $q'(s)=q(s)\cdot 
\prod_{i=0}^{r-1}\left(s+(i+n+\deg(g))/k\right)$. In particular,
\[
\left[\frakm^r\subseteq \fraka_{q(s)}\right]\Longrightarrow 
\left[b_Q(s)\divides
q(s)\cdot \prod_{i=0}^{r-1}\left(s+\frac{i+n}{k}\right)\right].
\] 
\proof
Let $m$ be a monomial of
degree $r-1$, so $x_img\in R_n[s]\cdot \fraka_{q(s)}$. Then 
\begin{eqnarray*}
\sum_{i=1}^n\del_i\action(x_imgQ^s)&=&
mg(\del_1x_1+\ldots +\del_nx_n)\action Q^s+\deg(mg)mgQ^s\\
&=&
mg(ks+n+\deg(mg))Q^s. 
\end{eqnarray*}
As $\deg(m)=r-1$, $\left(s+\frac{r-1+n+\deg(g)}{k}\right)mg\in \fraka_{q(s)}$.
By decreasing induction on $\deg(m)$,
\[
\prod_{i=1}^r\left(s+\frac{n+\deg(g)+r-i}{k}\right)
g\in
R_n[s]\cdot \fraka_{q(s)}.
\]

The final claim follows from the definition of $b_Q(s)$.
\qed
\end{lem}
\begin{rem}
  \label{rem-ex-nonarrangement}  
Suppose
that $f$ is {\em $w$-quasi-homogeneous}, i.e.,
  there are nonnegative numbers $w_1,\ldots,w_n$ such that with
  $\xi=\sum_{i=1}^n 
  w_ix_i\del_i$ one has $f=\xi\bullet(f)$ and hence $\xi-s\in J(f^s)$.  
If  $\frakn\subseteq\fraka_{q(s)}$ is a
  $w$-homogeneous 
  $\frakm$-primary ideal then one can show in the same manner that
$b_f(s)$ divides the
  product of 
  $q(s)$ and 
the minimal polynomial of $\xi$ on $R_n/\frakn$ evaluated
  at $-s-\sum_{i=1}^nw_i$. 
For example, $f=x^3+y^3+z^2w$ is (1/3,1/3,1/3,1/3)-homogeneous. 
One has $\fraka_{(s+1)}=\ideal{x^2,y^2,z^2,zw}$, which is of dimension
1, corresponding to the line of singularities 
$(0,0,0,w)$. One can see that the trick
of Lemma \ref{lem-euler} can be used to show that
$\fraka_{(s+1)(s+7/3)}=\ideal{x^2,xyz,y^2,z^2,zw}$ since $xyz$
is in the socle of $R_n/\fraka_{(s+1)}$. Going one step further,
$\fraka_{(s+1)(s+7/3)(s+2)}=\ideal{x^2,xz,y^2,yz,z^2,zw}$ and then $z$
can be obtained in
$\fraka_{(s+1)(s+7/3)(s+2)(s+5/3)}=\ideal{x^2,y^2,z}$. The new factors
are always equal to $s+\sum_{i=1}^4 (1/3)$ plus the degree of the new
element in $\fraka$. 

Now, however, nothing is in the socle and our
procedure stops. On the other hand, $f$ is also
$(1/3,1/3,1/2,0)$-homogeneous and this can be used to show that
\begin{eqnarray*}
\fraka_{(s+1)(s+7/3)(s+2)(s+5/3)(s+11/6)}&=&\ideal{x,y,z},\\
\fraka_{(s+1)(s+7/3)(s+2)(s+5/3)(s+11/6)(s+7/6)}&=&R_n.
\end{eqnarray*}
In fact,
$b_f(s)=(s+1)(s+7/3)(s+2)(s+5/3)(s+11/6)(s+7/6)$ and one can see again
how the factors of $b_f(s)$ enlarge (if taken in the right order) the
ideal $\fraka$, by either saturating, or dropping dimension.
\end{rem}

The trick for bounding $b_f(s)$ 
is therefore to find $q(s)$ such that $\fraka_{q(s)}$ is
zero-dimensional, and then to get a good estimate on the exponent $r$ of
Lemma \ref{lem-euler} if $g=1$. The importance of the relation $k\cdot
Q-\sum_{i=1}^nx_i\del_i$ in the annihilator of $Q^s$ for homogeneous
$Q$ of degree $k$ justifies 
\begin{df}
The {\em Euler operator} is $\euler=x_1\del_1+\ldots+x_n\del_n$. 
\end{df}

\subsection{Arrangements in the plane}\hfill

One has the following folklore result:
\begin{prop}
\mylabel{prop-in-plane}
Let $\{a_i\}_{3\le i\le k}$ be $k-2$ pairwise distinct nonzero numbers. Then
the Bernstein-Sato polynomial of $Q=xy(x+a_3y)\cdots(x+a_{k}y)$ divides
\[
(s+1)\prod_{i=0}^{2k-4}\left(s+\frac{i+2}{k}\right).
\] 
\proof
Consider the partial derivatives  $Q_x$ and $Q_y$ of $Q$ and 
the homogeneous forms
$x^iy^jQ_x$ and $x^iy^jQ_y$ where $i+j=k$. We claim that these
$2(k+1)$ forms of degree $(k+1)+k$ are linearly independent (and hence
that $\ideal{Q_x,Q_y}$ contains all monomials of degree at least $2k+1$). 

To see this, let $M=\{m_{a,b}\}_{0\le a,b\le 2k+3}$ 
be the matrix whose $(a,b)$-coefficient is the
coefficient of $x^{2k-3-b}y^b$ in $x^{k-2-a}y^aQ_y$ if $a\le k-2$, and
the coefficient of $x^{2k-3-b}y^b$ in $x^{2k-3-a}y^{a-k+1}Q_x$ if $a> k-2$.
%
The determinant of $M$ is the resultant of
$Q_x(1,y/x)$ and $Q_y(1,y/x)$. These cannot have a common root since
$\ideal{Q_x(x,y),Q_y(x,y)}$ is $\ideal{x,y}$-primary.
Hence  $M$ is of full rank and
$\frakm^k\ideal{Q_x,Q_y}=\frakm^{2k+1}$. 
Since $Q_x,Q_y\in\fraka_{(s+1)}$, $\frakm^{2k+1}\subseteq
\fraka_{(s+1)}$. Lemma \ref{lem-euler} implies the claim.
\qed
\end{prop}
Of course, a central  arrangement $Q$ of lines in the plane is 
an isolated singularity.
The interesting question was therefore the precise
determination of $\delta_Q$. 

\subsection{Estimates in dimension  $n>2$}\hfill

For the remainder of this section, $Q$ is a generic central
arrangement $Q=\prod_{i=1}^k H_i$. 
In order to estimate $b_Q(s)$ for $n>2,k>n+1$ we will consider 
a mix of the two main ideas for $n=2$. 
Namely,  we had $\frakm^{2k+1}\subseteq \fraka_{(s+1)}$. The point is
that admitting $(s+1)$ as a factor of $b_Q(s)$ allowed to capture
(set-theoretically) the singular locus of the arrangement.
This in 
conjunction with Lemma \ref{lem-euler} gave a bound for the
Bernstein-Sato polynomial. 

The plan is to devise a mechanism that starts with $\ideal{Q}\subseteq
\fraka_{1}$ and uses iterated multiplication with $(s+1)$ to enlarge
$\fraka_{q(s)}$. Progress is measured by the dimension of (the variety
of) 
$\fraka_{q(s)}$. This approach works well for generic arrangements,
while for non-generic arrangements or other singularities better
tricks seem to be needed.

It is crucial to understand the difference between the Jacobian
ideal of $Q$ and the ideal generated by
all $(n-1)$-fold products of distinct elements in $\calA$, and more
generally the difference between 
the ideal of the Jacobian ideal of the variety defined by all
$(r+1)$-fold products of distinct elements of $\calA$
and the ideal of all $r$-fold products of distinct elements
of $\calA$. 

\begin{df}
If $\calA=\{H_1,\ldots,H_k\}$ is a list of linear homogeneous
polynomials and $\alpha\in\NN^k$ we say that $\prod_{i=1}^k
{H_i}^{\alpha_i}$ is an {\em $\calA$-monomial}. If each $\alpha_i$ is either 0
or 1, we call the $\calA$-monomial {\em squarefree}. 
\end{df}

\begin{df}
We define polynomials $\Delta_{J,I,N}(Q)$ for 
a given central arrangement $Q=H_1\cdots H_k$. 
To this end let $N=\{\lambda_1,\ldots,\lambda_n\}\subseteq \{1,\ldots,k\}$ 
be a set of indices
serving as a
coordinate system. Let $v_{\lambda_1},\ldots,v_{\lambda_n}$ be $n$
appropriate $\CC$-linear 
combinations of $\del_1,\ldots,\del_n$ such that
$v_{\lambda_i}\action(H_{\lambda_j})=\delta_{i,j}$. 

Let $I\subseteq \{1,\ldots, k\}$ with  $|I|\geq k-n+1$. Set $\check I=
\{1,\ldots,k\}\setminus (I\cup N)$, 
$\hat I=I\cap N$, $H_I=\prod_{i\in I}H_i$. Observe
that $|\hat I|= |I|-k+n+|\check I|$.  

Let $\rho_N(I):=|\check I|+1\le |\hat I|$ and
pick $J\subseteq \hat I$ with $|J|=\rho_N(I)$. 
We define $\Delta_{J,I,N}(Q)$ to be
the $\rho_N(I)\times \rho_N(I)$-determinant and linear combination of
square-free $\calA$-monomials of degree $|I|-1$
\begin{eqnarray}
\label{eqn-def-Delta}
\Delta_{J,I,N}(Q)&=&\det\left(
\begin{array}{cccc}
v_{j_1}\action(H_{\check{\iota}_1})&\cdots&
v_{j_1}\action(H_{\check{\iota}_{|\check I|}}) &v_{j_1}\action(H_I)\\
\vdots&&\vdots&\vdots\\
v_{j_{|J|}}\action(H_{\check{\iota}_1})&\cdots&
v_{j_{|J|}}\action(H_{\check{\iota}_{|\check I|}}) &v_{j_{|J|}}\action(H_I)
\end{array}\right)
\end{eqnarray}
where $\check I=\{\check{\iota}_1,\ldots,\check\iota_{|\check
  I|}\}$ and $J=\{j_1,\ldots,j_{|J|}\}$. 
If $\check I$ is empty, $\Delta_{J,I,N}(Q)$ is just
  $\nu_{j_1}(H_I)$. 
We emphasize that $\Delta_{J,I,N}(Q)$ is defined only if $|I|>k-n$.
For a given $I$, let $\Delta_I(Q)$ be the set of all
$\Delta_{J,I,N}(Q)$, varying over all possible $N$, and for each $N$
over all $J$ satisfying  $J\subseteq\hat I$ and  $|J|=\rho_N(I)$.

Finally, put for $r>k-n$
\[
\Delta_r(Q)=\ideal{\Delta_{J,I,N}(Q):
|I|=r}+\ideal{H_I:|I|=r}
\]
and for all $r$
\[
\Sigma_{r}(Q)=\ideal{H_I : |I|=r}. 
\]
\end{df}
\begin{rem}
The ideal $\Sigma_{r-1}(Q)$ describes set-theoretically the locus where
simultaneously $k-r+2$ of the $H_i$ vanish (for the case $r<k-n+2$ see
see 
Lemma \ref{lem-k-n+1}), while $\Delta_r(Q)$ is the
Jacobian ideal of the variety to $\Sigma_{r}(Q)$. It is clear that
\[
\Sigma_{r-1}(Q)\,\,\supseteq\,\, 
\Delta_r(Q)\,\,=\,\,
\ideal{\{\Delta_{J,I,N}(Q):|I|=r\}}+\Sigma_{r}(Q)\,\,\supseteq\,\,
\Sigma_{r}(Q).  
\]
\end{rem}

The following is easily checked (since $Q$ is generic):
\begin{lem}
\mylabel{lem-k-n+1}
If $r\le k-n+1$ then $\Sigma_{r}(Q)=\frakm^r=\Delta_r(Q)$.\qed
\end{lem}

We can describe the ``difference'' of 
$\Delta_r(Q)$ and $\Sigma_{r-1}(Q)$ as follows:

\begin{prop}
\mylabel{prop-Delta-reduction}
Let $k\geq r\geq k-n+1$.
Then 
\[
\ann_{R_n}\left(\frac{\Sigma_{r-1}(Q)}{\Delta_r(Q)}\right)\supseteq
\frakm^{k-n}.
\] 
\proof
First let $k=n$ so that 
$n\geq r\geq 1$. In this case $\calA$-monomials and monomials are the same
concepts. Then $\Sigma_{r-1}(Q)$ is the
ideal of all squarefree ($\calA$-)monomials of degree $r-1$,  and
$\Delta_r(Q)$ is the ideal of all squarefree ($\calA$-)monomials of 
degree $r$
as well as all partial derivatives of these monomials. 
Clearly then  in this case
$\Delta_r(Q)=\Sigma_{r-1}(Q)$.

We shall prove the claim by induction on $k-n$ and we assume now that
$k>n$. 
Let $H_I$ be a squarefree
 $\calA$-monomial of degree $r$. We must show that $mH_I/H_i
\in \Delta_r(Q)$ for all $i\in I$ and all $m\in \frakm^{k-n}$. 

Pick $N\subseteq \{1,\ldots,k\}$ with $|N|=n$ and 
write $m=\sum_{j\in N} m_jH_j$, $m_j\in\frakm^{k-n-1}$. 
Consider the summands $m_jH_jH_I/H_i$ in $mH_I/H_i$. If $j=i$ or if 
$j\not \in I$, then certainly $H_jH_I/H_i\in \Delta_r(Q)$. Thus
we are reduced to showing that if $i\not =j\in I$ then
$m_jH_jH_I/H_i\in\Delta_r(Q)$. 

Note that if $k\geq r\geq
k-n+1$ then $k-1\geq r-1\geq k-1-n+1$.
Since $H_I/H_j$ is a squarefree 
$\calA\setminus \{H_j\}$-monomial of degree $(r-1)$ 
we may use the induction hypothesis on the
arrangement to $Q/H_j$ with $k-1\geq n$ factors. 
Hence for $i\not=j\in I$ there
 are $q_{j}\in 
\Delta_{r-1}(Q/H_j)$ such that $m_jH_I/H_iH_j=q_j$. Then
$m_jH_jH_I/H_i={H_j}^2q_j$ so that it suffices to show that  
\[
\left[q_j\in
\Delta_{r-1}(Q/H_j)\right] \Longrightarrow
\left[{H_j}^2q_j\in\Delta_r(Q)\right]. 
\]

It is sufficient to check this for $q_j$ being equal to
one of the two types of generators for $\Delta_{r-1}(Q/H_j)$, namely 
$H_{I'}$ and the 
determinants $\Delta_{J',I',N'}(Q/H_j)$, where as usual $J',I',N'\subseteq
\{1,\ldots,k\}\setminus \{j\}$, $|N'|=n$ and $|I'|=r-1$. If 
$q_j=H_{I'}$ then
 ${H_j}^2q_j=H_jH_{I'\cup\{j\}}\in\Delta_r(Q)$. So assume that
$q_j=\Delta_{J',I',N'}(Q/H_j)$.

Multiplication of $\Delta_{J',I'
,N'}(Q/H_j)$ by ${H_j}^2$ can be achieved by multiplying the last column of the
defining matrix (\ref{eqn-def-Delta}) 
 of $\Delta_{J',I',N'}(Q/H_j)$ by ${H_j}^2$. 
Let in that context 
$j_t\in N'$ and $v_{j_t}$ be the corresponding derivation relative to
$N'$. 
Then
\[
{H_j}^2v'_{j_t}\action(H_{I'})=
H_jv'_{j_t}\action(H_{I'\cup\{j\}})-H_{I'\cup\{j\}}v'_{j_t}\action(H_j).
\]
Thus, ${H_j}^2\Delta_{J',I',N'}(Q/H_j)=H_j\Delta_{J',I',N'}(Q)$ modulo
$\ideal{H_{I'\cup\{j\}}}$. As $H_{I'\cup\{j\}}\in \Delta_r(Q)$, 
${H_j}^2\Delta_{r-1}(Q/H_j)\subseteq \Delta_r(Q)$.
The proposition follows hence by induction.
\qed
\end{prop}

Recall that $\Delta_I(Q)$ is the collection of  all
$\Delta_{J,I,N}(Q)$ for fixed $I$. We now relate the ideals
$\Sigma_{r-1}(Q)$ 
and
$\Delta_r(Q)$ to ideals $\fraka_{q(s)}$ and give hence the latter
ideals geometric meaning.
\begin{lem}
\mylabel{lem-Deltaprime-reduction}
Fix integers $r\geq k-n+2$, and $t$. 
Suppose $\frakm^tH_I\subseteq \fraka_{q(s)}$ for some $I$ with $|I|=r$. 
Then $\frakm^{t+1}\Delta_I(Q)\subseteq
\fraka_{(s+1)q(s)}$. In particular,
\[
\left[\frakm^t\Sigma_{r}(Q)\subseteq \fraka_{q(s)}\right]\Longrightarrow
\left[\frakm^{t+1}\Delta_r(Q)\subseteq \fraka_{(s+1)q(s)}\right].
\] 
\proof
%
Pick a specific $\Delta_{J,I,N}(Q)$ and a monomial $m$ of degree $t$. 
In particular, this means that a coordinate system
$H_{\lambda_1},\ldots,H_{\lambda_n}$ and derivations 
$v_{\lambda_1},\ldots,v_{\lambda_n}$ have been chosen. 
Consider the effect of $v_j$ on $mH_IQ^s$ for $j\in J\,(\subseteq I\cap N)$:
\begin{eqnarray*}
\frac{v_j\action\left(mH_IQ^s\right)}{Q^s}&=&\frac{1}{Q}
  \left(v_j\action(m)H_IQ+
   mQv_j\action(H_I)+
   smH_Iv_j\action(Q)
   \right)\\ 
&=&v_j\action(m)H_I+
   (s+1)mv_j\action(H_I)+
   \sum_{i\in\{1,\ldots,k\}\setminus I}
   smH_I\frac{v_j\action(H_i)}{H_i}
\end{eqnarray*} 
The sum has only poles of order one.  These poles occur exactly 
along all hyperplanes in
$\check I$ since $v_j\action(H_i)=0$ if $i\not =j, i\in N$. (Note
that $j\in J\subseteq I$  is not index of a summand.) 
The $|\check I|+1$ distinct elements of $J$ give rise to that
many expressions of the type shown.
Hence there is a nontrivial $\CC$-linear combination of the
$v_j\action(mH_IQ^s)$ without poles; by construction this linear
combination is in
$\fraka_{q(s)}$. 
It is easy to see that the
desired expression results in
$(s+1)m\Delta_{J,I,N}(Q)+V(m)H_I$ where $V(m)$ is a
linear combination in the $v_j\action(m)$. 
As $x_iv_j\action(m)H_I \in\frakm^tH_I\subseteq \fraka_{q(s)}$, 
$x_im\Delta_{J,I,N}(Q)\in \fraka_{(s+1)q(s)}$ for all $i,J,N$ and so
$\frakm^{t+1}\Delta_I(Q)\subseteq \fraka_{(s+1)q(s)}$. 

To prove the final assertion, note that $\Delta_r(Q)$ is generated by
all 
$\Delta_I(Q)$, $|I|=r$ and all $H_I$, $|I|=r$. One then only needs to
observe that all $H_I$ with $|I|=r$ are already in $\Sigma_{r}(Q)$. 
\qed
\end{lem}

One can now conclude alternately 
from \ref{prop-Delta-reduction} and
\ref{lem-Deltaprime-reduction} 
that
\begin{eqnarray*}
\Sigma_k(Q)&\subseteq&\fraka_1,\\
\frakm \Delta_k(Q)&\subseteq&\fraka_{(s+1)},\\
\frakm^{k-n+1} \Sigma_{k-1}(Q)&\subseteq&\fraka_{(s+1)},\\
&\vdots&\\
\frakm^{(k-n+1)(n-2)+1}\Delta_{k-n+2}(Q)&\subseteq&\fraka_{(s+1)^{n-1}},\\
\frakm^{(k-n+1)(n-1)}\Sigma_{k-n+1}(Q)&\subseteq&\fraka_{(s+1)^{n-1}},
\end{eqnarray*}
and since $\Sigma_{k-n+1}(Q)=\frakm^{k-n+1}$,
$\frakm^{(k-n+1)n}\subseteq\fraka_{(s+1)^{n-1}}$. 
It is very intriguing how in the above sequence of containments an
  extra factor of $(s+1)$ in $q(s)$ 
allows each time to reduce the dimension of
  $\fraka_{q(s)}$ and in fact to enlarge $\fraka_{q(s)}$ to an ideal
  with radical equal to the singular locus of
  $\fraka_{q(s)}$. One might compare this to the example in Remark
  \ref{rem-ex-nonarrangement}. 

The remainder of this section is devoted to decreasing substantially the
exponent of $\frakm$ in the final row of the display above.

\begin{prop}
\mylabel{prop-Deltaprime-cup-m}
For all $r\in \NN$ with $k-n+1\le r\le k+1$, 
\[
\frakm^{2k-n-1}\cap\Sigma_{r-1}(Q)\subseteq \fraka_{(s+1)^{k-r+1}}.
\] 
\proof
We shall proceed by decreasing induction on $r$. We know that
\[
\frakm^{2k-n-1}\cap\Sigma_{k}(Q)\subseteq \Sigma_{k}(Q)=
\ideal{Q}\subseteq\fraka_{(s+1)^0}.
\]
 Assume then that $k-n+1\le r\le k$ and that
$\frakm^{2k-n-1}\cap\Sigma_{r}(Q)\subseteq \fraka_{(s+1)^{k-r}}$. 
Since $\Sigma_{r}(Q)\subseteq \frakm$ is homogeneous of degree $r$, 
this implies that
\[
\frakm^{2k-n-1-r}\cdot \Sigma_{r}(Q)\subseteq
\fraka_{(s+1)^{k-r}}.
\]

We
need to show that
$\frakm^{2k-n-1}\cap\Sigma_{r-1}(Q)\subseteq\fraka_{(s+1)^{k-r+1}}$ in
order to get the induction going. For this, we consider $\Delta_r(Q)$. 
Let $\Delta$ be a generator of 
$\Delta_r(Q)$. Either $\Delta=H_I$ and  $|I|=r$, in which case
$\Delta\in \Sigma_{r}(Q)$.
Or, $\Delta=\Delta_{J,I,N}(Q)$ with $|I|=r$. In that case, Lemma
\ref{lem-Deltaprime-reduction}  
together with 
$\frakm^{2k-n-1-r}\cdot\Sigma_{r}(Q)\subseteq
\fraka_{(s+1)^{k-r}}$
implies that 
$\frakm^{2k-n-r}\cdot \Delta\subseteq 
\fraka_{(s+1)^{k-r+1}}$.
Therefore our hypotheses imply that 
\[
\frakm^{2k-n-1}\cap\Delta_r(Q)\subseteq\fraka_{(s+1)^{k-r+1}}.
\] 
But then,
\begin{eqnarray*}
\frakm^{2k-n-1}\cap\Sigma_{r-1}(Q)&=&
\frakm^{2k-n-1-(r-1)}\Sigma_{r-1}(Q)\\
&& \text{ ($\Sigma_{r-1}(Q)$ is
homogeneously generated in degree $r-1$)}\\
&=&\frakm^{k-r}\frakm^{k-n}\Sigma_{r-1}(Q)\\
&\subseteq &\frakm^{k-r}(\Delta_r(Q)\cap\frakm^{k-n+r-1})\\
&& \text{ (by Proposition \ref{prop-Delta-reduction})}\\
&\subseteq&\frakm^{2k-n-1}\cap\Delta_r(Q)\\
&\subseteq&\fraka_{(s+1)^{k-r+1}}.
\end{eqnarray*}
\qed
\end{prop}
This proposition says that sufficiently high degree parts of the ideal
defining the higher iterated singular loci of $\calA$ are 
contained in certain $\fraka_{q(s)}$. It gives quite directly a bound
for the 
Bernstein-Sato polynomial:
\begin{thm}
\mylabel{thm-upper-bound}
The Bernstein-Sato polynomial of the central generic arrangement 
$Q=H_1\cdots H_k$ divides
\begin{eqnarray}
\label{eqn-upper-bound}
(s+1)^{n-1}\prod_{i=0}^{2k-n-2}\left(s+\frac{i+n}{k}\right).
\end{eqnarray}
\proof
The previous proposition shows (with $r=k-n+2$) that
$\frakm^{2k-n-1}\cap \Sigma_{k-n+1}\subseteq
\fraka_{(s+1)^{n-1}}$. By Lemma \ref{lem-k-n+1},
$\Sigma_{k-n+1}(Q)=\frakm^{k-n+1}$. Thus,
$\frakm^{2k-n-1}\subseteq\fraka_{(s+1)^{n-1}}$. We conclude now as in
Lemma \ref{lem-euler}.\qed
\end{thm}
In the next two sections we show that this estimate is in essence the
correct answer.

\section{Remarks on a  conjecture by Orlik and Randell}
Let $Q:\CC^n\to \CC$ be a homogeneous 
polynomial map, denote by $X_\alpha$ the preimage
$Q^{-1}(\alpha)$ for $\alpha\in \CC\setminus \{0\}$ and let  $X$ be the fiber
over zero. As $Q$ is homogeneous 
the $X_\alpha$ 
are all isomorphic and smooth. 
 Let $\tilde \CC^\times$ be the universal
cover of $\CC^\times=\CC\setminus \{0\}$, and $\tilde X$ the fiber product
of $\tilde \CC^\times$ and $\CC^n\setminus X$ over $\CC^\times$. Then 
$(\alpha,x)\to
(\alpha+2\pi,x)$ is a diffeomorphism of $\tilde X$ and therefore induces an
isomorphism  
$\mu$ on the cohomology $H^*(X_\alpha,\CC)$, the Picard-Lefschetz
monodromy \cite{Brieskorn,Deligne,Hamm}. 
If in addition $X$ has an isolated
singularity then $X_\alpha$ is homotopy equivalent to a
bouquet of $(n-1)$-spheres \cite{Milnor} and so the only (reduced) 
cohomology of
the fiber is in degree $n-1$. The roots of the minimal
polynomial $a_\mu(s)$ of $\mu$
are in that case obtained from the roots of the
Bernstein-Sato polynomial of $Q$ by $\lambda\to e^{2\pi i\lambda}$ 
\cite{Malgrange}. The multiplicities remain mysterious, however. 
If $X$ is not an isolated singularity, the $X_\alpha$ have cohomology in
degrees other than $n-1$ and the monodromy acts on all these
cohomology groups. The monodromy is then not so nicely related to the
Bernstein-Sato polynomial and not well understood.

\subsection{The conjecture}\hfill

The natural projection $R_n\onto
R_n/\ideal{Q-\alpha}$ induces a map of 
differentials $\Omega\to\Omega_\alpha$
which in turn induces a surjective map of de Rham complexes
$\pi:\Omega^\bullet\onto \Omega^\bullet_\alpha$ where
$\Omega_\alpha$ 
are the $\CC$-linear differentials on $R_n/\ideal{Q-\alpha}$ and
$\Omega^\bullet_\alpha$ is the de Rham complex on $X_\alpha$. 
It is an interesting and open
question to determine explicit formul\ae\ for generators of the
cohomology of $\Omega^\bullet_\alpha$, i.e.\ forms on $\CC^n$ that restrict
to generators of $H^i(\Omega_\alpha^\bullet)$, $i\le n-1$. 
If $X$
has an isolated singularity then the Jacobian ideal $\frakA$ is
Artinian, the dimension of the vector space $R_n/\frakA$ equals
$\dim_\CC(H^{n-1}_{\DR}(X_\alpha,\CC))$, and
 the elements of $R_n/\frakA$ can be
identified with the 
classes in $H^{n-1}_{\DR}(X_\alpha,\CC)$. Namely, $\bar g\in R_n/\frakA$ 
corresponds to
$g\omega$ where 
\begin{eqnarray}
\label{eqn-omega}
\omega=\sum_{i=1}^n(-1)^{i-1}x_i\,dx_1\wedge\ldots\wedge
\widehat{dx_i}\wedge\ldots\wedge dx_n
\end{eqnarray}
and the hat indicates omission.
\bigskip

For the remainder of this section let 
$Q$ be a reduced polynomial describing a central generic
arrangement, $Q=H_1\cdots H_k$.
We let as before $\calA=\{H_1,\ldots,H_k\}$.
In \cite{Orlik-Randell} (Proposition 3.9) 
it is proved that every cohomology class in
$H^{n-1}(\Omega^\bullet_\alpha)$ is of the form $\bar{\pi(g\omega)}$  
for some $g\in R_n$, and that
\[
\dim_\CC(H^{n-1}(\Omega^\bullet_\alpha))= {k-2\choose n-2}+k{k-2\choose
n-1}.
\] 
The authors make a conjecture which states roughly that $g$ may be
chosen to be homogeneous  and that Milnor fibers
of central generic arrangements have a cohomology description similar
to the isolated singularity case. 

By $(R_n)_r$ we denote the homogeneous elements in $R_n$ of degree
$r$. The following vector space is central to the ideas of Orlik and
Randell.
\begin{df}
We denote by $\mu$ a subset of $\calA$ of cardinality $n-1$. We write
then $J_\mu(a)$ with $a\in R_n$ for the Jacobian determinant associated
to $H_{\mu_1},\ldots,H_{\mu_{n-1}},a$. We also  denote by $Q_\mu$ the
product of all $H_i$ with $i\not \in \mu$, its degree is hence $k-n+1$.
In our previous notation, $Q_\mu $ was $H_I$ with
$I=\calA\setminus\mu$. 

With these notations, let $E$ be the vector space in $R_n$ generated by
all elements of the form 
\begin{eqnarray}
\label{eqn-E-gens}
\deg(a)a\,J_\mu(Q_\mu)-kQ_\mu \,J_\mu(a),
\end{eqnarray}
varying over all homogeneous $a\in R_n$. It is not an $R_n$-ideal.
\end{df}

\begin{conj}[Orlik-Randell, \cite{Orlik-Randell}]
\mylabel{conj-o-r}
Consider the fiber $X_1=\var(Q-1)$. There is a finite dimensional
homogeneous vector space $U\subset R_n$ such that
\begin{enumerate}
\item $R_n= E\oplus (\CC[Q]\otimes U)$; 
\item the map $U\to H^{n-1}(X_1,\CC)$ given by $g\to
\bar{\pi(g\omega)} $ is an isomorphism, and
$\Omega^{n-1}_\alpha=\pi(U\omega)\oplus d\Omega^{n-2}_\alpha$.
\item The dimensions $u_r$ of $U_r$, the graded pieces of $U$ of degree $r$,
are as follows:
\[
u_r=\left\{\begin{array}{ll}
{r+n-1\choose n-1}&\text{for $0\le r\le k-n$},\\
{k-2\choose n-1}&\text{for $k-n+1\le r\le k-1$},\\
{k-2\choose n-1}-{r-k+n-1\choose n-1}&\text{for $k\le r\le 2k-n-2$}.
\end{array}
\right.
\]
\end{enumerate}
\end{conj}

In this section we will prove that if $k$ does not divide $r-k+n$ then
the dimension of $(R_n/E)_r$ 
is bounded by ${k-1\choose n-1}$ and that
strict inequality holds if additionally 
 $r>k$. In the next section we will see
that $(R_n)_r=E_r+\ideal{Q}_r$ for $r\geq 2k-n-1$. 
This will imply that $(R_n/(E+\ideal{Q-1}))_r$ is
nonzero exactly if  $0\le r\le 2k-n-2$, and that for $k-n-1\le r\le k$ its
dimension is exactly as the conjecture by Orlik and Randell  
predicts. 

It is worth pointing out that the vector space $E$ is too small if $Q$
is an arrangement that is not generic. For example, with
$Q=xyz(x+y)(x+z)$ as in Example \ref{ex-non-generic} one obtains that
the dimension of $(R_n/E+\ideal{Q})_r$ is 2 whenever
$r$ is at least $5$. 

\subsection{Generators for $U$}\hfill

We now consider the question of finding generators for $U$. 
By Lemma 
\ref{lem-k-n+1}, 
$(R_n)_{k-n+1}$ is generated by the set of 
all $Q_\mu$ as a vector space.
Then $(R_n)_r$ is for $r>k-n+1$ generated by
$\frakm^{r-k+n-1}\cdot \Sigma_{k-n+1}(Q)$. 
We claim that we may pick vector space generators $G=\{g_i\}$ for
$(R_n)_r$, $r>k-n+1$,
such that 
\begin{enumerate}
\item[a)] each $g_i$ is an $\calA$-monomial,
\item[b)] each $g_i$ is a
multiple of some $Q_\mu$.
\end{enumerate}
To see this, observe that  $
(R_n)_r=\left(\frakm^{r-k+n-1}\right)_{r-k+n-1}\cdot
\left(\Sigma_{k-n+1}(Q)\right)_{k-n+1}$. Since $\calA$ is essential, 
 Lemma \ref{lem-k-n+1} completes the argument.
We call an element of $R_n$ satisfying
these two conditions a {\em standard product}.

We shall now prove that there are no more than ${k-2\choose n-1}$ 
standard products
necessary to generate $(R_n/(E+\ideal{Q-1}))_r$. For $k-n+1\le r< k$
this is exactly the number 
stipulated by 
Conjecture \ref{conj-o-r}. We will do this by showing that the
relations in $E$ may 
be used to eliminate the majority of all summands in a typical element
of $(R_n)_r/E_r$. In order to do this, we need to study the nature of the
relations in $E$. To get started, note that 
\[
\left[H_j\in\calA\setminus\mu\right]\Longrightarrow \left[J_\mu(H_j)\not
  =0\right]. 
\]
We will now show that every generator (\ref{eqn-E-gens}) of $E$
induces  a syzygy
between $k-n+1$ squarefree $\calA$-monomials of degree $k-n$.
\begin{lem}
\mylabel{lem-relation-structure}
Let $a\in (R_n)_r$ 
be an $\calA$-monomial of positive
degree $r$ such that $k\dividesnot r$, and pick $n-1$ distinct factors
$\mu$ of $Q$.  
Consider the corresponding element 
\begin{eqnarray}
\label{eqn-relation-structure}
\deg(a)aJ_\mu(Q_\mu)-kQ_\mu J_\mu(a) 
\end{eqnarray}
of $E$.
In this expression
(using the product rule for computing the Jacobian) the first term
contributes $k-n+1$ summands of the form $\deg(a)a\frac{Q_\mu}{H_i}
J_\mu(H_i)$ where $H_i$ runs through the factors of
$Q_\mu$. Similarly the second term contributes $\deg(a)$ summands of
the form $kQ_\mu\frac{a}{a_i}J_\mu(a_i)$ with $a_i$ running through the
factors of $a$. We claim that all 
nonzero summands in the latter set (apart from
constant factors)  appear as  nonzero summands in the former set. Moreover,
for each summand that is nonzero on both sides the coefficients are
different. 
\proof
There are two main cases: $a_i\in\mu$ and $a_i\not\in \mu$. If
$a_i\in\mu$ then $J_\mu(a_i)$ is a determinant with a repeated column,
and hence the summand $Q_\mu\frac{a}{a_i}J_\mu(a_i)$ is zero. On the
other hand, 
$H_i\not\in\mu$ gives a summand
 $\deg(a)a\frac{Q_\mu}{H_i}J_\mu(H_i)\not= 0$.
So the left term in (\ref{eqn-relation-structure}) gives $k-n+1$ nonzero
$\calA$-monomials with nonzero coefficients.
If
$a_i\not \in\mu$, then $a_i=H_j$ (say), and
$\frac{Q_\mu}{H_j}a=Q_\mu\frac{a}{a_i}$. Let $t$ be the multiplicity
of $H_j$ in $a$, $a=a'\cdot {H_j}^t$. In (\ref{eqn-relation-structure}) the
first term contributes $\deg(a)Q_\mu\frac{a}{H_j}J_\mu(H_j)$ while
the second yields $t$ times $-kQ_\mu\frac{a}{H_j}J_\mu(H_j)$ by the
product rule. So the total number of copies of
$\frac{aQ_\mu}{H_j}J_\mu(H_j)$ in (\ref{eqn-relation-structure}) is
$\deg(a)-kt$.

As $k$ is not a divisor of $\deg(a)=r$, 
 each generator of $E_r$ gives rise to a relation between exactly
$k-n+1$ 
of
our generators of $(R_n)_r$, corresponding to the divisors of $Q_\mu$. 
\qed 
\end{lem}
\begin{rem}
Suppose that in a linear combination of $\calA$-monomials the previous
lemma is used to eliminate $Q_\mu\frac{a}{H_i}J_\mu(H_i)$. Then the
replacing $\calA$-monomials are of the form
$Q_\mu\frac{a}{H_j}J_\mu(H_j)$ where $H_j\not\in\mu$. 
\end{rem}
We now show how to use Lemma \ref{lem-relation-structure} to limit the
dimension of $(R_n)_r/E_r$.

\begin{prop}
\mylabel{prop-nice-gens}
Let $r\in\NN$, $k-n+1\le r$, and $k\dividesnot (r-k+n)$. 
The (cosets of) 
$\calA$-monomials of the form
\begin{eqnarray}
\label{eqn-nice-gens}
H_{i_1}\cdots H_{i_{k-n-1}}H_{k-1}{H_k}^{r-k+n}, \quad i_1<\ldots
<i_{k-n-1}<k-1 
\end{eqnarray}
span $(R_n/E)_r$ and therefore generate
$\left(H^{n-1}_{\DR}(Q^{-1}(1),\CC)\right)_r$.  
\proof
Let $P\in (R_n)_r$ be a standard product.
We prove that it may be replaced by a linear combination of
$\calA$-monomials
of the stipulated form. Here are three ways of modifying a linear
combination of $\calA$-monomials modulo $E$:
\begin{enumerate}
\item If $P$ uses $l>k-n+1$ distinct factors of $\calA$ we can
write  $P=P'Q_\mu$ for a suitable $\mu$ and we can assume that $H_k\in
\mu$. That means that $H_k\dividesnot Q_\mu$ and the multiplicity of
$H_k$ in $P'$ is of course at most $r-k+n-1$.
Let
$i_0=\min\{i:H_i\not \in\mu\}$
and $\mu'=\mu\cup \{H_{i_0}\}\setminus
\{H_{k}\}$, so $Q_{\mu'}=H_kQ_\mu/H_{i_0}$. Consider the element of $E$ 
given by $(r-k+n)P'H_{i_0}J_{\mu'}(Q_{\mu'})-kQ_{\mu'}
J_{\mu'}(P'H_{i_0})$. It is a linear dependence modulo
$E$ between $P'H_{i_0}Q_{\mu'}/H_k=P$ on one side and terms of the
form $P'H_{i_0}Q_{\mu'}/H_i=P'H_kQ_\mu/H_i$ 
for $H_k\not =H_i\in\calA\setminus \mu'$
on the other, 
with no
coefficient equal to zero. It follows that  $P=P'Q_\mu$ may,
modulo $E$, be replaced by a linear combination of standard products 
with a higher power of $H_k$ in each of them than in
$P$ and $l$ or $l-1$   distinct factors. Note that each replacing
$\calA$-monomial has multiplicity of $H_k$ at most $r-k+n$.

\item Suppose now that $P$ has exactly $k-n+1$ distinct factors, but
that $H_k$ is not one of them. Let $Q_\mu$ be the product of all
distinct factors of $P$, and set $P=P'Q_\mu$. Let
$i_0=\min\{i:H_i\not\in\mu\}$ 
and set 
$\mu'=\mu\cup \{H_{i_0}\}\setminus 
\{H_{k}\}$. 
The relation $(r-k+n)H_{i_0}P'J_{\mu'}(Q_{\mu'})-kQ_{\mu'}
J_{\mu'}(P'H_{i_0})$ allows to replace $P$ by a linear combination of
standard products with $k-n+1$ or $k-n+2$ distinct factors (depending
on the multiplicity of $H_{i_0}$ in $P'$)  such that $H_k$
divides each of the new standard products.
\item Now assume that $P$ is a standard product with exactly $k-n+1$ distinct
  factors and assume furthermore
that $H_k$ divides $P$ with multiplicity $l<r-k+n$. Let $\mu$ be
such that $Q_\mu$ divides $P$.
Since the arrangement is
generic, the $n-1$ elements of $\mu$, 
together with $H_k$, span the maximal ideal and
thus if $i_0=\min\{i:{H_i}^2\divides P\}$ then one  factor $H_{i_0}$ of
  $P$ may be
replaced by an appropriate linear combination in $H_k$ and the
  elements of $\mu$. 
This creates a linear combination of $(n-1)$ standard products 
with $k-n+2$
distinct factors in each summand where $H_k$ has
  multiplicity $l$, and one $\calA$-monomial with $k-n+1$ factors
  where the $H_k$-degree is $l+1$. 
\end{enumerate}
Starting with any standard product of degree $r$, using these steps in
appropriate order will produce a linear combination of 
standard products with exactly $k-n+1$ factors and 
multiplicity $r-k+n$ in $H_k$. This is
because after every
execution of Step 1 and 2, the degree in $H_k$ goes up, and 
after each execution of Step 3 we may do Step 1 at least once on the
$n-1$ standard products with $k-n+1$ factors.

Now let $P=H_{i_1}\cdots H_{i_{k-n}}{H_k}^{r-k+n}$ 
with $i_1<i_2\cdots<i_{k-n}<k-1$. Let $\mu$ be such that
$Q_\mu=H_{i_1}\cdots H_{i_{k-n}}H_{k-1}$, in particular $H_k\in\mu$. 
Then
\[
E\ni
(r-k+n){H_k}^{r-k+n}J_{\mu}(Q_{\mu})-kQ_{\mu}J_{\mu}({H_k}^{r-k+n})=
  (r-k+n){H_k}^{r-k+n}J_{\mu}(Q_{\mu})
\]
allows to replace $P$ by a sum of $\calA$-monomials 
 each of which has $k-n+1$
distinct $\calA$-factors, and each of which is divisible by
$H_{k-1}{H_k}^{r-k+n}$ (note that the only term that might fail to
have $H_{k-1}$ in it disappears because $J_\mu({H_k}^{r-k+n})=0$ as
$H_k\in\mu$). Thus, modulo $E$, $P$ is equivalent to a linear
combination of $\calA$-monomials of type (\ref{eqn-nice-gens}).

The condition $k\dividesnot (r-k+n)$ is needed because otherwise Lemma
\ref{lem-relation-structure} does not work.
\qed
\end{prop}

\begin{rem}
Note that there are exactly ${k-2\choose k-n-1}={k-2\choose n-1}$
$\calA$-monomials of type (\ref{eqn-nice-gens}). It follows that 
$\dim(R_n/E)_r\le {k-2\choose n-1}$ unless $k$ divides $r-k+n$. 
Also, if $r=k-n$ the conjecture says that ${k-2\choose n-1}$
generators for $(R_n/E)_r$ are not enough. So in a sense this is an
optimal estimate.
In the
following section we will see that $(R_n)_r=E_r$ along $Q^{-1}(1)$ for
$r>2k-n-2$. We finish this section with a lemma that will be used in
the next section to prove that $(R_n)_r\not = E_r$ for $k\le r\le 2k-n-2$.
\end{rem}

\begin{lem}
\mylabel{lem-nontrivial-relation}
If $r\geq k$ and $k\dividesnot (r-k+n)$ then 
$\dim_\CC(R_n/E+\ideal{Q})_r\le {k-2\choose n-1}-1$.
\proof
The proof of Proposition \ref{prop-nice-gens} contains a procedure to
turn $Q{H_k}^{r-k}$ into a sum of $\calA$-monomials of the form
(\ref{eqn-nice-gens}). One may do so using only Step 1 of that proof. In
fact, if $P=H_1\cdots H_{k-n-1}H_{i_1}\cdots
H_{i_{j}}H_{k-1}{H_k}^{r-k+n-j}$ for $k-n-1<i_1<\cdots <i_j<k-1$, then the
relation (\ref{eqn-relation-structure}) induced by $Q_\mu=H_1\cdots
H_{k-n-1}H_{i_1}H_k$ and $a=P/Q_\mu$ allows to replace $P$ by a sum of
$\calA$-monomials
 each
of which is divisible by $H_{k-1}{H_k}^{r-k+n+1-j}$, each of which has
only $H_k$ as repeated factor, and precisely one of which is a nonzero
multiple of
$H_1\cdots H_{k-n-1}$. Therefore rewriting $Q{H_k}^{r-k}$ only using
Step 1 (and only $Q_\mu=H_1\cdots H_{k-n-1}\cdot H_{i_1}H_k$ with
$k-n\le i_1<k$) gives a
relation modulo $E$ 
between the products of (\ref{eqn-nice-gens}) where the
coefficient for $H_1\cdots H_{k-n-1}H_{k-1}{H_k}^{r-k+n}$ is
nonzero. Hence in particular, 
$(R_n/E+\ideal{Q})_r$ has dimension at most ${k-2 \choose n-1}-1$. 
\qed
\end{lem}
We have shown that
filtering $R_n/(E+\ideal{Q-1})$ by degree, the $r$-th graded piece
has dimension at most ${k-2\choose n-1}-1$ unless $k$ divides
$r-k+n$. Moreover, $(R_n/E)_r=(R_n)_r$ for 
$r\le k-n$. 

\section{Integration, Restriction and Bernstein-Sato polynomials}

If $Q$ is radical and describes a generic arrangement then we have
seen that
\begin{itemize}
\item $b_Q(s)$ is a divisor of
$(s+1)^{n-1}\prod_{i=0}^{2k-n-2}(s+\frac{i+n}{k})$; 
\item $\dim_\CC(R_n/E+\ideal{Q})_r\le {k-2\choose n-1}$ if
  $k\dividesnot (r-k+n)$; 
\item the inequality of the previous item is strict if  in addition $r<k-n$ or 
  $r\geq k$.
\end{itemize}
We now prove, among other things, 
 that each homogeneous
generator of the top cohomology group 
$H^{n-1}_{\DR}(Q^{-1}(1),\CC)$ of the Milnor fiber 
gives for all homogeneous
 polynomials $Q$ rise to a root of $b_Q(s)$
just as it is the case for homogeneous isolated
singularities. 

\subsection{Restriction and Integration}\hfill

A central part of this section is occupied by effective methods for
$D$-modules. In fact, we shall use in an abstract way algorithms that
were pioneered by Oaku \cite{Oa2} and since have become the
centerpiece of algorithmic $D$-module theory.

We shall first explain some basic facts about restriction and
integration functors. 
Much more detailed explanations may be found in
\cite{Oa2,O-T2,O-T1} and \cite{W-cdrc}. In particular, we only consider the
situation of $n+1$ variables $x_1,\ldots,x_n,t$ and explain
restriction to $t=0$ and integration along $\del_1,\ldots,\del_n$.

\begin{df} 
Let $\tilde\Omega_t=D_{x,t}/t\cdot D_{x,t}$ and
$\Omega_\del=D_{x,t}/\{\del_1,\ldots,\del_n\}\cdot D_{x,t}$. 

The {\em restriction of the $D_{x,t}$-complex
  $A^\bullet$ to the 
subspace $t=0$} is the complex $\rho_t(A^\bullet)=\tilde\Omega_t\otimes
_{ D_{x,t} }^LA^\bullet$ 
considered as a complex in the category
of $D_x$-modules.

The {\em integration of $A^\bullet$ along
$\del_1,\ldots,\del_n $} is the
complex $\DR(A^\bullet)=\Omega_\del\otimes 
_{ D_{x,t} }^LA^\bullet$ considered as a complex in
the
category of $D_t$-modules. 
\end{df}
In the sequel we describe tools that may be used to compute
restriction and integration.

\begin{dfs}

\mylabel{V-definitions}
On the ring $ D_{x,t} $  the {\em $V_t$-filtration} 
$F_t^l( D_{x,t} )$ is the $\CC$-linear span of all
operators $x^\alpha\del^\beta t^a\del_t^b$ for 
which $a+l\geq b$.
More
generally, on a free $ D_{x,t} $-module 
$A=\bigoplus_{j=1}^r  D_{x,t} \cdot e_j$ we set 
\[
F^l_t(A[\frakm])=
\sum_{j=1}^r
F^{l-\frakm(j)}_t( D_{x,t} )\cdot e_j, 
\]
where $\frakm$ is any element of $\ZZ^r$ called  the {\em shift
vector}. A shift vector is tied to a fixed set of generators.
The {\em $V_t$-degree} 
of
an operator $P\in A[\frakm]$ is the
smallest $l=V_t\ddeg (P[\frakm])$ such that 
$P\in F^l_t(A[\frakm])$.

If $M$ is a quotient of the free $ D_{x,t} $-module $A=\bigoplus_1^r
D_{x,t} \cdot 
e_j$,  
$M=A/I$, 
we define the $V_t$-filtration on
$M$ by $F^l_t(M[\frakm])=F^l_t(A[\frakm])+I$ and
for submodules $N$ of $A$ by
intersection: $F^l_t(N[\frakm])=F^l_t(A[\frakm])\cap N$.



\end{dfs}
\begin{dfs}
\mylabel{graded-dfs}
A complex of free $ D_{x,t} $-modules 
\[
\cdots\to
A^{i-1}\stackrel{\phi^{i-1}}\longrightarrow
A^i\stackrel{\phi^i}\longrightarrow 
A^{i+1}\to\cdots
\]
 is said to be {\em $V_t$-strict } with respect to
the shift vectors $\{\frakm_i\}$ if
\[
\phi^i\left(F^l_t(A^i[\frakm_i])\right)\subseteq F^l_t(A^{i+1}[\frakm_{i+1}])
\] 
and also 
\[
\im(\phi^{i-1})\cap
F^l_t(A^i[\frakm_i])=\im(\phi^{i-1}|_{F^l_t(A^{i-1}[\frakm_{i-1}])})
\]
for all $i,l$.

Set $\theta=t\del_t$, the Euler operator for $t$.  
A $ D_{x,t} $-module $M[\frakm]=A[\frakm]/I$ 
is called {\em specializable to $t=0$} if
there is a polynomial $b(s)$ in a single variable such that
\begin{eqnarray}
\label{eqn-spec}
b(\theta+l)\cdot F^l_t(M[\frakm]) &\subseteq& F^{l-1}_t(M[\frakm])
\end{eqnarray}
for all $l$ (cf.\
\cite{K2,O-T2}). Holonomic modules are specializable.
Introducing 
\[
\gr^l_t(M[\frakm])=(F^l_t(M[\frakm]))\//\/(F^{l-1}_t(M[\frakm])),
\] 
this can be
written as 
\[
b(\theta+l)\cdot \gr^l_t(M[\frakm])=0.
\]
The monic polynomial $b(\theta)$ of least degree satisfying an
equation of the type (\ref{eqn-spec}) is called the {\em
$b$-function for restriction of $M[\frakm]$ to $t=0$}.
\end{dfs}
By \cite{O-T1}
(Proposition 3.8) and \cite{W-cdrc} every complex 
admits a $V$-strict resolution.
In the theorems to follow, the meaning of filtration on restriction
and integration complexes is as in \cite{W-cdrc}, Definition 5.6.

\begin{thm}[\cite{Oa2,O-T1,W-cdrc}]
\mylabel{alg-restriction}
Let $(A^\bullet[\frakm_\bullet],\delta^\bullet)$ 
be a $V_t$-strict complex of free
$ D_{x,t}$-modules with holonomic cohomology. The restriction
$\rho_t(A^\bullet[\frakm_\bullet])$ of $A^\bullet[\frakm_\bullet]$
to $t=0$ can
be computed as follows:
\begin{enumerate}
\item Compute the $b$-function $b_{A^\bullet[\frakm_\bullet]}(s)$ for
restriction of $A^\bullet[\frakm_\bullet]$ to $t=0$.
\item Find an integer $l_1$ with $\left[b_{A^\bullet[\frakm_\bullet]}(l)=0, 
l\in\ZZ\right]\Rightarrow
\left[l\le l_1\right]$. 
\item $\rho_t(A^\bullet[\frakm_\bullet])$ is 
quasi-isomorphic to the complex 
\begin{eqnarray}
&\cdots\to F^{l_1}_t(\tilde\Omega_t\otimes_{ D_{x,t} }A^i[\frakm_i])
\to
F^{l_1}_t(\tilde\Omega_t\otimes_{ D_{x,t} }A^{i+1}[\frakm_{i+1}])
\to
\cdots
&
\label{eqn-restricted-complex}
\end{eqnarray}
\end{enumerate}
This is a complex of free finitely generated 
$D_x$-modules and a representative of
$\rho_t(A^\bullet[\frakm_\bullet])$. Moreover, if a cohomology class
in $\rho_t(A^\bullet[\frakm_\bullet])$ has $V_t$-degree $d$ then $d$
is a zero of $b_{A^\bullet[\frakm_\bullet]}(s)$.
\qed
\end{thm}
In order to compute the integration along $\del_1,\ldots,\del_n$ 
one defines a filtration by 
\[
F^l_\del(D_{x,t})=\{x^\alpha\del^\beta t^a\del_t^b: |\alpha|\le
|\beta|+l\}.
\]
With $\tilde \euler=-\del_1x_1-\ldots-\del_nx_n$ the 
{\em $b$-function for
integration} of the module $M$ is the 
least degree monic polynomial $\tilde b(s)$ such that
\[
\tilde b(\tilde \euler+l)\cdot F^l_\del(M)\subseteq F^{l-1}_\del(M).
\]
Then the integration complex $\DR(M)$ of $M$ is quasi-isomorphic
to 
\begin{eqnarray}
&\cdots\to{\tilde F}^{l_1}_\del(\Omega_\del\otimes_{ D_{x,t} }A^i[\frakm_i])
\to{\tilde F}^{l_1}_\del(\Omega_\del\otimes_{ D_{x,t} }A^{i+1}[\frakm_{i+1}])
\to
\cdots
&
\label{integrated-complex}
\end{eqnarray}
where $A^\bullet[\frakm_\bullet]$ is a $V_\del$-strict resolution of
$M$, and $l_1$ is the largest integral root of $\tilde b(s)$. Again,
cohomology generators have $V_\del$-degree equal to a root of $\tilde
b(s)$. 

\subsection{Bernstein-Sato polynomial and the relative de Rham
complex}\hfill

Following Malgrange \cite{Malgrange2}, 
we consider for $f\in R_n$ 
the symbol $f^s$ as generating a $D_{x,t}$-module contained
in the free $R_n[f^{-1},s]$-module $R_n[f^{-1},s]f^s$ via 
\[
t\action \frac{g(s)}{f^j}f^s=\frac{g(s+1)}{f^{j-1}}f^s,\qquad
\del_t\action \frac{g(s)}{f^j}f^s=\frac{-sg(s-1)}{f^{j+1}}f^s.
\]
Then the left ideal $J_{n+1}(f)=\ideal{t-f,\{\del_i+
\del_i\action(f)\del_t\}_{i=1}^n}\subseteq D_{x,t}$
is easily seen to consist of operators that
annihilate $f^s$. Moreover, $-\del_tt$ acts as multiplication by $s$.
Since $J_{n+1}(f)$ is maximal, it
actually contains all annihilators of $f^s$. It turns out, that
$J_{n+1}(f)$ describes the $\calD$-module direct image of $R_n$ under the
embedding $x\to (x,f(x))$:
\begin{lem}
For all $f\in R_n$, 
\[
D_{x,t}/J_{n+1}(f)\cong H^1_{t-f}(R_{x,t}),
\]
generated by $\frac{1}{t-f}$.
\proof
Consider $\tau=\frac{1}{t-f}\in R_{x,t}[(t-f)^{-1}]$. It is obviously
annihilated by $\{\del_i+
\del_i\action(f)\del_t\}_{i=1}^n$. Moreover, $\frac{t-f}{t-f}\in
R_{x,t}$ so that $(t-f)(\tau\mod R_{x,t})=
0\in H^1_{t-f}(R_{x,t})$. Hence
$J_{n+1}(f)$ annihilates the coset of $\tau$ in $H^1_{t-f}(R_{x,t})$.  

The polynomial $t-f$ is free of singularities and so its
Bernstein-Sato polynomial is $s+1$. Hence $\tau$ generates
$R_{x,t}[(t-f)^{-1}]$. Therefore 
the coset of $\tau$ generates the local cohomology module. Since this
module is nonzero, the coset of $\tau$ cannot be zero. Hence its
annihilator cannot be $D_{x,t}$. As $J_{n+1}(f)$ is maximal 
we are done.\qed
\end{lem}

We now connect the ideas of Malgrange with algorithmic methods
pioneered by Oaku, and Takayama, to show that
the ideal $J_{n+1}(f)$ is intimately
connected with the Bernstein-Sato polynomial of $f$:

\begin{thm}
\mylabel{thm-bernstein-integration}
Let $Q$ be a homogeneous polynomial of degree $k>0$ with 
Bernstein-Sato polynomial
$b_Q(s)$. Then $b_Q((-s-n)/k)$
is a multiple of the $b$-function for integration of
$J_{n+1}(Q)=\ideal{t-Q,\{\del_i+\del_i\action(Q)\del_t\}_{i=1}^n}$ along
$\del_1,\ldots,\del_n$.  
\proof
It is well-known, that $J_{n+1}(Q)\cap D_x[s]=\ann_{D_x[s]}(Q^s)$. 
Hence in particular, $J_{n+1}(Q)$ contains
$\euler-ks=\euler+k\del_tt$. 

To be a Bernstein polynomial means that $b_Q(s)\in
J_{n+1}(Q)\cap D_x[s]+D_x[s]\cdot Q$.  Write $b_Q(s)=j+P(s)Q$ with 
$j\in J_{n+1}(Q)\cap D_x[s]$, $P(s)\in
D_x[s]$. 

The ideal $J_{n+1}(Q)$ is $(1,k)$-homogeneous if we set $\deg(x_i)=1$,
$\deg(\del_i)=-1$, $\deg(t)=k$, $\deg(\del_t)=-k$. Since $b_Q(s)$ is
$(1,k)$-homogeneous of degree 0, 
we may assume that $j$ (and hence $P(s)Q$) is also
$(1,k)$-homogeneous
of degree
$0$. Writing $P(s)=\sum_{i=0}^l P_is^i$ with $P_i\in D_x$, we see
that each $P_i$ is of $(1,k)$-degree $-k$. 
This implies, as $P_i\in D_x$, that $P_i\in
F^{-k}_{\del}(D_x)$. Note that as $t-Q\in J_{n+1}(Q)$, $b_Q(s)=P(s)t$
modulo $J_{n+1}(Q)$. So $P(s)t=P(-\del_tt)t\in F^{-k}_{\del}(D_{x,t})$
and $b_Q(-\del_tt)\in J_{n+1}(Q)+F^{-k}_\del(D_{x,t})$.

Also, $b_Q(-\del_t t)$ is modulo $J_{n+1}(Q)$ equivalent to $b_Q((-\tilde
\euler-n)/k)$ because  $\euler+k\del_tt\in J_{n+1}(Q)$. Thus
\[
b_Q\left(\frac{-\tilde \euler-n}{k}\right)\in 
J_{n+1}(Q)+F^{-k}_{\del}(D_{x,t}),
\]
proving that $b_Q(-(s+n)/k)$ is a multiple of the $b$-function for
integration of $D_{x,t}/J_{n+1}(Q)$ 
along $\del_1,\ldots,\del_n$.\qed
\end{thm}
Combining Theorem \ref{thm-bernstein-integration} with Theorem
\ref{alg-restriction} and its integration counterpart, one obtains
\begin{cor}
\label{cor-bernstein-integration}
The only possible $V_\del$-degrees  for the generators of the
cohomology of $\DR(D_{x,t}/J_{n+1}(Q))$ 
are those specified by the roots of
$b_Q(-(s+n)/k)$.\qed
\end{cor}

\subsection{Restriction to the fiber}\hfill

Let $Q$ be a homogeneous polynomial of positive degree. We now consider
the effect of restriction to $t-1$ of the relative de Rham complex
$\DR(D_{x,t}/J_{n+1}(Q))$.
This is computed as the cohomology of
the tensor product over $D_t$ of $\DR(D_{x,t}/J_{n+1}(Q))$ with
$(D_{x,t}\stackrel{(t-1)\cdot}{\longrightarrow} D_{x,t})$. We shall
concentrate on the highest cohomology group. It equals
$D_{x,t}/(J_{n+1}(Q)+\{\del_1,\ldots,\del_n,t-1\}D_{x,t})$. 


\begin{thm}
\label{thm-U-spanned}
The quotient 
\begin{eqnarray}
\label{eqn-def-U}
U&:=&D_{x,t}/(J_{n+1}(Q)+\{\del_1,\ldots,\del_n,t-1\}D_{x,t})
\end{eqnarray}
 is spanned by 
polynomials $g\in R_n$. One may choose a $\CC$-basis for $U$ in
such a way that 
\begin{itemize}
\item  all basis elements are in $R_n$ and homogeneous;
\item  no basis element may be replaced by an element of smaller
degree (homogeneous 
or not).
\end{itemize}
We call such a basis a {\em homogeneous degree minimal basis}.

The degree of any element $g$ of a degree minimal basis satisfies 
\[
b_Q(-(\deg(g)+n)/k)=0
\]
and then the usual degree of $g$ is the $V_\del$-degree of the class of
$g$. 
\proof
Clearly $U$ is spanned by the cosets of $R_n[\del_t]$. Let $g\in R_n$
be homogeneous of degree $d$. In $U$ we have $gt^a\del_t^b=gt^b\del_t^b=
g\prod_{j=0}^{b-1}(t\del_t-j)$ for all
$a,b\in\NN$. Now observe that $\euler+k\del_tt\in J_{n+1}(Q)$ implies that
in $U$
\begin{eqnarray*}
0=\del_t^bg(\euler+k\del_tt)&=&\del_t^b(\euler-d+k\del_tt)g \\
&=&\del_t^b(-n-d+k\del_tt)g\\
&=&\left((-n-d+k(b+1))\del_t^b+k\del_t^{b+1}\right)g.
\end{eqnarray*}
By induction this shows that in $U$ 
\begin{eqnarray}
\label{eqn-tdt-acts}
\del_t^bg=t^b\del_t^bg=\prod_{i=1}^b\left(\frac{n+d}{k}-i\right)g.
\end{eqnarray}
Hence $U$ is spanned by the cosets of $R_n$. As $\del_t^bg$ and
$\prod_{i=1}^b\left(\frac{n+d}{k}-i\right)g$ have the same
$V_\del$-degree, minimal $V_\del\ddeg$
representatives for all $u\in U$ can be chosen within $R_n$.

Now let $u'\in R_n$ be homogeneous and let $0\not =u\in D_{x,t}$ 
be a $V_\del$-degree minimal
representative of the class of $u'$ in $U$. 
By the previous paragraph, without affecting 
$V_\del$-degree, $u$ can be assumed to be in $R_n$.
Then obviously the $V_\del$-degree agrees with the usual degree. 

Therefore by definition
$\deg(u)=V_\del\ddeg(u)\le V_\del\ddeg(u')=\deg(u')$. Hence 
for any $u'\in U$ we have 
\[
\min\{\deg(u):u\in R_n, u=u'\in U\}=\min\{V_\del\ddeg(u):u'=u\in U\}
\]
and the equality can be realized by one and the same element $u\in
R_n$ on both 
sides. If this $u$ is nonzero in $U$, then clearly $u\in
H^n(\DR(D_{x,t}/J_{n+1}(Q)))$ is also nonzero since this module
surjects onto $U$. The $V_\del$-degree of
$u$ within $H^n(\DR(D_{x,t}/J_{n+1}(Q)))$ cannot be smaller than the
$V_\del$-degree of $u$ in the bigger coset when considered in $U$, and
hence is just the usual degree of $u$. 
By Corollary \ref{cor-bernstein-integration}, 
$b_Q(-(\deg(u)+n)/k)=0$.
This implies that $U$ is finite dimensional. Hence any $\CC$-basis for
$U$ may be turned into a degree minimal one. 

It remains to show that the basis can be picked in a homogeneous
minimal way.
Note that
$J_{n+1}(Q)+\{\del_1,\ldots,\del_n,t-1\}D_{x,t}$ is
$\ZZ/\ideal{k}$-graded (by $x_i\mapsto 1,\del_i\mapsto -1$ and
$t,\del_t\mapsto 0$); so $U$ is $\ZZ/\ideal{k}$-graded and $U$ has a
degree minimal $\ZZ/\ideal{k}$-graded basis.  If  
$u$ is in a $\ZZ/\ideal{k}$-graded minimal degree basis but 
not homogeneous, the degrees of its 
graded components only
differ by multiples of $k$. Write $u=u_a+u_{a+1}+\ldots+u_b$
with $a,b\in \NN$ and $u_j$ the component of $u$ in degree
$jk$. Then since $t-Q$ is in $J_{n+1}(f)$, we have in $U$ the equality 
$u=\sum_{j=a}^bu_jQ^{b-j}$. The right hand side is
homogeneous and both of usual and of $V_\del$-degree $\deg(u)$. 
Hence $\ZZ/\ideal{k}$-graded 
minimal degree bases for $U$ can be changed into homogeneous
minimal degree bases without changing the occurring degrees (all of
which we proved to be roots of $b_Q(-(s+n)/k)$).
\qed
\end{thm}
\begin{rem}
\mylabel{rem-restriction-at-1}
An important hidden ingredient of the above theorem is the fact that
the 
$b$-function for restriction to $t-1$ of both $D_{x,t}/J_{n+1}(f)$ and
$H^n \DR(D_{x,t}/J_{n+1}(f)) $ is
$(t-1)\del_t$ whenever $f$ is $w$-homogeneous. 
Namely,
if  $f=\sum_{i=1}^n
w_ix_i\del_i\action{f}$ then with $\xi=\sum_{i=1}^n
w_ix_i\del_i$ we have $(\del_tt+\xi)\action f^s=0$.
Consider then the equation
\[
(t-1)\del_t=\underbrace{(t-1)(\del_tt+\xi)}_A-
\underbrace{(t-1)(\del_t(t-1)+\xi)}_B.  
\]
Obviously, $A\in J_{n+1}(f)$ and $B\in F^{-1}_{t-1}(D_{x,t})\cap
F^{-1}_{\del,t-1}(D_{x,t})$.  These are the required conditions to be a
$b$-function for restriction to $t-1$ of $D_{x,t}/J_{n+1}(f)$
respectively $H^n(\DR(D_{x,t}/J_{n+1}(f)))$.
\end{rem}
%

\begin{cor}
\mylabel{cor-U-spanned}
Let $Q=H_1\cdots H_k$ define a central generic arrangement. Then $U$
has a homogeneous basis of polynomials of degree at most $2k-n-2$.
\proof
By Theorem \ref{thm-upper-bound},  $b_Q(s)$ has its zero locus
inside $\{-n/k,\ldots,(-2k+2)/k\}$. Then by Theorem
\ref{thm-U-spanned}  the degrees of a minimal degree basis for $U$
are bounded above by $2k-n-2$.
\qed
\end{cor}

\subsection{De Rham cohomology from $\calD$-module operations}\hfill

For any $f$, 
the complex  $\DR(D_{x,t}/J_{n+1}(f))$ carries the de Rham cohomology of the
fibers of the map $\CC^n\ni x\to f(x)\in\CC$, since it is the result of
applying the de Rham functor to the composition of maps $x\to(x,f(x))$
and $(x,y)\to (y)$ (see \cite{Deligne}). 
The de Rham functor for the embedding
corresponds to the functor that takes the $D_x$-module $M$ to the
$D_{x,t}$-complex $M\otimes_{D_x}
(D_{x,t}\stackrel{\cdot(f-t)}{\longrightarrow}D_{x,t})$, while the
projection corresponds to the formation of the Koszul complex induced
by left multiplication by $\del_1,\ldots,\del_n$. 
The cohomology of the fiber $Q^{-1}(1)$ is obtained as the restriction to
$t-1$. 

With the shifts in cohomological degree, $U=D_{x,t}/
(J_{n+1}(Q)+
\ideal{\del_1,\ldots,\del_n,t-1}\cdot D_{x,t})$ thus encodes the top de Rham
cohomology of $Q^{-1}(1)$. For homogeneous $Q$ 
the correspondence between these two spaces
is at follows. Write $dX=dx_1\wedge\ldots \wedge dx_n$ and
$\widehat{dX_i}=dx_1\wedge\ldots\wedge\widehat{dx_i}\wedge\ldots\wedge
dx_n$ where the hat indicates omission. 
An element $g$ in $U$ determines the form
$g\,dX$ on $\CC^n$. Under the embedding $Q^{-1}(1)\into
\CC^n$, the form $g\, dX$ restricts ($\calD$-module
theoretically) to the form
$G$ which satisfies $dQ\wedge G=g\, dX$.   Let us
compute $G$. Since $G$ is an $(n-1)$-form, $G=\sum_{i=1}^n
g_i\,\widehat{dX_i}$. 
Thus, $dQ\wedge
G=\sum_{i=1}^n(-1)^i\del_i\action(Q)
g_i\, dX$. On the other hand, along $Q^{-1}(1)$,
$g\, dX=gQ\, dX=\frac{g}{k}\sum_{i=1}^nx_i\del_i\action(Q)
dX$.
 Thus by comparison, $kg_i=(-1)^ix_ig$.
With $\omega$ as in (\ref{eqn-omega}), 
the $(n-1)$-form on $Q^{-1}(1)$ 
encoded by $g\in U$ is $G=g\omega/k$.  We show now that all forms in
$H^{n-1}_{\DR}(Q^{-1}(1),\CC)$ are captured by $U$. 

\begin{lem}
\mylabel{lem-omega}
If $Q$ is homogeneous, then $H^{n-1}_{\DR}(Q^{-1}(1),\CC)$ is
generated by $R_n\cdot\omega$.
\proof This is trivial for $n=1$, so we assume that $n>1$. 
Consider the map $R_n\to H^{n-1}_{\DR}(Q^{-1}(1),\CC)$ given by $g\to
g\omega$. Suppose $g\omega=0$. Then 
\[
g\omega=(Q-1)h+d(G)+A\wedge dQ
\]
for $h=\sum (-1)^{i+1}h_i
\,\widehat{dX_i}\in\Omega^{n-1}$,  
$G=\sum g_{i,j}
dx_1\wedge\ldots\wedge\widehat{dx_i}\wedge\ldots\wedge \widehat{dx_j}
\wedge\ldots\wedge 
dx_n \in \Omega^{n-2}$, $A\in \Omega^{n-2}$. Multiply by
$dQ$ to get
\[
kQg\,dX=\left(\sum_{i=1}^n
(Q-1)\del_i\action(Q)h_i+\sum_{i,j=1}^n\del_i\action(Q)\del_j\action(g_{i,j})
-\del_j\action(Q)\del_i\action(g_{i,j})\right)dX
\]
in $\Omega^n=R_n\,dX$. Now look at this in
$U$. Note that $kQg=ktg=kg$ and
$(Q-1)\del_i\action(Q)h_i=(t-1)\del_i\action(Q)h_i=0$ in $U$. So (in
$U$)  
\begin{eqnarray}
\label{eqn-skew-stuff}
kg=\sum_{i,j=1}^n\del_i\action(Q)\del_j\action(g_{i,j})-
 \del_j\action(Q)\del_i\action(g_{i,j})
\end{eqnarray}
We would like this to be zero in $U$; in fact it will turn out to
vanish term by term.

We may assume that $g_{i,j}$ is homogeneous by looking at the graded
pieces $g$ of (\ref{eqn-skew-stuff}). So to simplify notation let $h$
be a homogeneous polynomial in $R_n$. In the remainder of this proof we
shall use a subscript to denote derivatives:
$h_i=\del_i\action(h)$. 
Then in $U$ we have
\begin{eqnarray*}
0=-\del_jth_i+\del_ith_j&=&
-th_{i,j}+th_{j,i}+(-th_i\del_j+th_j\del_i)\\
&=&th_iQ_j\del_t-th_jQ_i\del_t\\
&=&(h_iQ_j-h_jQ_i)t\del_t\\
&=&(h_iQ_j-h_jQ_i)\frac{n+\deg(h)-2}{k}\qquad\text{by (\ref{eqn-tdt-acts}).}
\end{eqnarray*}
If $\deg(h)>0$, this implies the vanishing of $h_iQ_j-h_jQ_i\in U$.
But if $\deg(h)=0$ there is nothing to  prove in the
first place. Therefore the sum (\ref{eqn-skew-stuff}) is zero.
 Hence  if $g\omega=0$ in
$H^{n-1}_{\DR}(Q^{-1}(1),\CC)$ then $kg=0$ in $U$. So $R_n\onto U$
factors as $R_n\onto R_n\omega\onto U=H^{n-1}_{\DR}(Q^{-1}(1),\CC)$. \qed
\end{lem}

Our considerations prove in view of Corollary \ref{cor-U-spanned}:
\begin{thm}
\mylabel{thm-deRham-gens}
Let $Q\in R_n$ be a homogeneous polynomial of degree $k$.
The de Rham cohomology group $H^{n-1}_{\DR}(Q^{-1}(1),\CC)$ is 
isomorphic to 
$U\cdot \omega$. There is a homogeneous basis for $U$ with degrees
bounded by 
\[
u_Q=\max\{i\in\ZZ:b_Q(-(i+n)/k)=0\}.
\]
If $Q$ defines a generic arrangement of hyperplanes, $u_Q\le 2k-n-2$.
\qed
\end{thm}

\subsection{Non-vanishing of  $H^{n-1}_{\DR}(Q^{-1}(1),\CC)$ and roots
  of $b_Q(s)$}\hfill 

We now establish the existence of a non-vanishing 
$g\in H^{n-1}_{\DR}(Q^{-1}(1),\CC)$  
in all degrees
$0\le\deg(g)\le 2k-n-2$ for generic central arrangements $Q$. This
will certify each root of (\ref{eqn-upper-bound}) as root of $b_Q(s)$.

The primitive 
$k$-th root $\zeta_k$ of unity acts on $\CC^n$ by $x_i\to \zeta_k x_i$. 
This fixes $Q$ and hence the ideal
$J_{n+1}(Q)$. Therefore it gives an automorphism of the de Rham complex
and hence the induced map on cohomology
separates $H^{n-1}_{\DR}(Q^{-1}(1),\CC)$ into eigenspaces,
$U=\bigoplus_{\bar i\in\ZZ/k\ZZ}M_{\bar i}$ which are
classified by their degree modulo $k$.

From
\cite{Orlik-Randell} we know the monodromy of $Q$.
In particular, $M_{\bar i}$ is
a ${k-2\choose n-1}$-dimensional vector
space unless $\bar i=\bar{k-n}$. 
Write $U_i$ for the elements in $U$ with (homogeneous) 
minimal degree representative
of degree precisely $i$.
Since 
elements of $U$ 
have degree
at most $2k-n-2$, we find that 
\begin{eqnarray*}
\dim(U_i)+\dim(U_{i+k})=\dim(U_i+U_{i+k})={k-2\choose n-1}&\text{ for
}&0\le i\le k-n-1;\\
U_i=M_{\bar i}\text{ and } \dim(U_i)={k-2\choose n-1}&\text{ for }&
             k-n<i\le k-1.
\end{eqnarray*}
Moreover, $U_{k-n}=M_{\bar {k-n}}=(R_n)_{k-n}$ of dimension
${k-1\choose n-1}$.

Since $R_n/(E+\ideal{Q-1})$ surjects onto $U$, 
Lemma \ref{lem-nontrivial-relation}  shows that neither $U_i$ nor $U_{i+k}$
is zero-dimensional for $0\le i\le k-n-1$. So one has 
\begin{thm}
For a generic hyperplane arrangement $Q$ the vector space 
\begin{eqnarray*}
\left(H^{n-1}_{\DR}(Q^{-1}(1),\CC)\right)_r\not = 0 \,\,\text{ for
}\,\,0\le r\le 
2k-n-2.
\end{eqnarray*}
It is zero for all other $i$.\qed
\end{thm}

One can now use the non-vanishing to certify roots of $b_Q(s)$ as such:
%
%

\begin{cor}
\mylabel{cor-result}
The Bernstein-Sato polynomial of a generic central arrangement
$Q=\prod_{H_i\in \calA}H_i$
of degree $k$ is 
\[
(s+1)^{r}\prod_{i=0}^{2k-n-2}\left(s+\frac{i+n}{k}\right)
\]
where $r=n-1$ or $r=n-2$.
\proof
By the previous theorem, $U_i\not =0$ for $0\le i\le 2k-n-2$. A
minimal degree basis for $U$ must therefore contain elements of all
these degrees. By the last part of Theorem \ref{thm-U-spanned},
$b_Q(s)$ is a multiple of
$\prod_{i=0}^{2k-n-2}\left(s+\frac{i+n}{k}\right)$. On the other hand,
Theorem \ref{thm-upper-bound} proves that $b_Q(s)$ divides the
displayed expression with $r=n-1$. This proves everything apart from
the multiplicity of $(s+1)$.

Let $\vec x\not=\vec 0$ 
be any point of the arrangement where precisely $n-1$ planes meet.
 The Bernstein-Sato polynomial of $Q$ is a multiple of the
local Bernstein-Sato polynomial at $\vec x$ (which is defined by the same
type of equation as $b_Q(s)$ but where $P(s)$ is in the localization
of $D_x[s]$ at the maximal ideal defining 
$\vec x$). Since the local Bernstein-Sato polynomial at
a normal crossing of $n-1$ smooth divisors is $(s+1)^{n-1}$, the theorem
follows. \qed
\end{cor}
We are quite certain, that the exponent $r$ in Corollary
\ref{cor-result} is $n-1$, but we do not know how to show that.
In fact, we believe that the elements $g\in R_n$ whose cosets in 
\[
(s+1)^{n-2}\prod_{i=0}^{2k-n-2}\left(s+\frac{i+n}{k}\right)
\cdot \frac{D_n[s]\action
  f^s}{D_n[s]\action f^{s+1}}
\]
are zero are precisely the elements of $\frakm^{k-n+1}$.

\begin{conj}
If $k\le r\le 2k-n-2$ we believe that 
the space $(R_n/(E+\ideal{Q-1}))_r$ is
spanned by the expressions in (\ref{eqn-nice-gens}) for which
$i_1<(n-1)+(r-k)$. If $k-n<r<k$, the expressions in Proposition
\ref{prop-nice-gens} are known to span $U$.  If $r\le k-n$ we believe
that $U_r=(R_n)_r$. 

This is in accordance with \cite{Orlik-Randell} as 
there are exactly as many such expressions
as Conjecture \ref{conj-o-r} predicts for 
the dimension of $(H^{n-1}_{\DR}(Q^{-1}(1),\CC))_r$.
\end{conj}

\begin{ex}
\mylabel{ex-non-generic}
Consider the non-generic arrangement given by
$Q=xyz(x+y)(x+z)$. With the $D$-module package \cite{M2D} of {\sl
  Macaulay 2} \cite{M2} one computes its Bernstein-Sato polynomial
as
\[
(s+1)(s+\frac{2}{3})(s+\frac{3}{3})(s+\frac{4}{3})
(s+\frac{3}{5})(s+\frac{4}{5})(s+\frac{5}{5})(s+\frac{6}{5})
(s+\frac{7}{5}).
\]
Therefore the $b$-function for integration of $J_{n+1}$ along
$\del_1,\ldots,\del_n$ is a divisor of 
\[
(s-2)(s-\frac{1}{3})(s-2)(s-\frac{11}{3})(s-0)(s-1)(s-2)(s-3)(s-4).
\]
This indicates that the degrees of the top cohomology of the Milnor fiber
$Q^{-1}(1)$ are at most 4. It also shows that in this case these
degrees do not suffice to determine the roots of $b_Q(s)$. In fact,
the degrees of no class in any $H^i_{\DR}(Q^{-1}(1),\CC)$ will explain
the roots $-2/3$ and $-4/3$ in $b_Q(s)$. 

However, consider  a point $P\not =0$ on 
the line 
$x=y=0$. This line is the intersection of three participating
hyperplanes, $x$, $y$ and $x+y$. In $P$ the variety of $Q$ has a
homogeneous structure as well, so the local Bernstein-Sato polynomial
of $Q$ at $P$ is a multiple of the minimal polynomial of the local
Euler operator on the cohomology of the Milnor fiber of $Q$ at
$P$. In fact, at $P$ the variety of $Q$ is a generic arrangement in
the plane, times the affine line. Without difficulty one verifies then
that the Milnor fiber has top cohomology in degrees 0, 1 and 2, and
that 
$b_{Q,P}(s)=(s+2/3)(s+1)^2(s+4/3)$. 

The global Bernstein-Sato polynomial of $Q$ is the least common
multiple of all local Bernstein-Sato polynomials $b_{Q,P}(s)$. Hence
$b_Q(s)$ must be a multiple of $(s+2/3)(s+1)^2(s+4/3)$ and so all
roots of $b_Q(s)$ come in one way or another from cohomology degrees
on Milnor fibers. This prompts
the following question:
\begin{prob}
Let $Q$ be a locally quasi-homogeneous polynomial in $R_n$ (as for
example a hyperplane arrangement).
Is it true that every root of $b_Q(s)$ arises through the action of an
Euler operator on the top de Rham cohomology of the Milnor fiber of
$Q$ at some point of the arrangement?
\end{prob}

This is of course true for isolated quasi-homogeneous
singularities. If $Q$ is an arrangement then
by the local-to-global principle one may restrict to central
arrangements. We have proved here that the question has an affirmative
answer for generic arrangements.

One more remark is in order.
The cohomology we have used to describe the Bernstein-Sato polynomial
is the one with coefficients in the constant sheaf $\CC$, which may
be viewed as 
the sheaf of solutions of the $D_n$-ideal
$\ideal{\del_1,\ldots,\del_n}$ describing $R_n$ on $\CC^n$. Relating
holonomic $D_n$-modules to locally constant sheaves on $\CC^n$ is the point of
view of the  Riemann-Hilbert correspondence, \cite{Boreletal}.
There are, however, other natural locally constant sheaves on
$\CC^n\setminus 
Q^{-1}(0)$ induced by $D_n$-modules 
than just the constant sheaf. For example, for every
$a\in\CC$ the $D_n$-ideal $\ann_{D_n}(f^a)$ induces such a sheaf as
the sheaf of its local solutions. For most exponents $a$ this is of
course a sheaf without global sections on $Q\not =0$, 
and more generally without any
cohomology. For suitable exponents, however, this is different; it is
sufficient to consider the case where
$a+\ZZ$ contains a root of the Bernstein-Sato polynomial.
Perhaps one can
characterize the Bernstein-Sato polynomial as the polynomial of
smallest degree such that $s=-\del_tt$ annihilates the $V_\del$-degree
of every cohomology class in $H^i(\Omega\otimes_{D_n}^L\int_\iota\calP)$ for
every $D_n$-module $\calP$ defining a locally constant system on
$\CC^n\setminus Q^{-1}(0)\stackrel{\iota}{\hookrightarrow}
\CC^{n+1}$. Another possibility is given by the cyclic covers introduced by
Cohen and Orlik \cite{Cohen-Orlik}.
\end{ex}


\section{Miscellaneous Results}

In this section we collect some results and conjectures 
concerning the structure of
the module $D_n\action Q^s$ associated to central arrangements.

\subsection{Arbitrary arrangements}\hfill

We begin with a fact pointed out to us by A.\ Leykin.
\begin{thm}[Leykin]
\mylabel{thm-leykin}
The only integral root of the Bernstein-Sato polynomial of any
arrangement $\calA$  is $-1$.
\proof
By Lemma \ref{lem-broot}
it will be sufficient to show that if $Q=\prod_{H_i\in \calA}H_i$ then
$R_n[Q^{-1}]$ is generated by $1/Q$ since this implies that $D_n\action
(Q^{-1})=D_n\action(Q^{-r})$ for all $r\in\NN$. 

Since the Bernstein-Sato polynomial is the least common multiple of
the local Bernstein-Sato polynomials, we may assume that $\calA$ is
central.
We may also assume that
$\calA\subseteq \CC^n$
is not contained in a linear subspace of $\CC^n$.

The claim
is true for a normal crossing arrangement.
We proceed by induction on the difference $k-n>0$ where
$k=\deg(Q)$. Since the local cohomology module $H^k_\frakm(R_n)$ vanishes,
$R_n[Q^{-1}]=\sum_{i=1}^k R_n[(Q/H_i)^{-1}]$. Moreover, by induction
$R_n[(Q/H_i)^{-1}]$ is
generated by $H_i/Q$ as $D_n$-module. Since obviously $H_i/Q$ is in
the $D_n$-module generated by $1/Q$, the theorem follows.
\qed
\end{thm}

\begin{rem}
Note that the same argument proves the following. Let
$g_1,\ldots,g_k\in R_n$ and set $G=\prod_{i=1}^kg_i$. If $R_n[(G/g_i)^{-1}]$ is
generated by $(g_i/G)^m$ for $i=1,\ldots,k$, and
$H^k_{\ideal{g_1,\ldots,g_k}}(R_n)=0$ then $R_n[G^{-1}]$ is generated by
$1/G^m$. That is to say, if the smallest integral root of
$b_{G/g_i}(s)$ is at least $-m$, and if
$H^k_{\ideal{g_1,\ldots,g_k}}(R_n)=0$ then the smallest integral root
of $b_{G}(s)$ is at least $-m$. By Grothendieck's vanishing theorem
this last condition is always satisfied if $k>n$. 
\end{rem}

We now give some combinatorial results on the localization module
$R_n[Q^{-1}]$. 
The following really is a general fact about finite length modules.
\begin{prop}
\mylabel{prop-in-ex}
Let $M=\sum_{i=1}^kM_i$ be a holonomic $D_n$-module. Then the
holonomic length  satisfies
\[
\ell(M)=\sum_{i=1}^k(-1)^{i+1}\sum_{|I|=i}\ell(M_I)\]
where $M_I=\bigcap
_{j\in I}M_j$.
\proof
$\ell$ is additive in short exact sequences. Hence
$\ell(M)=\ell(M_1)+\ell(M/M_1)$. 
In order to start the induction, one needs to look at the case $k=2$
which is the second isomorphism theorem.

Also, by induction,
\begin{eqnarray*}
\ell(M)-\ell(M_1)=\ell(M/M_1)
&=&\ell(\sum_{i>1}(M_j+M_1)/M_1)\\
&=&\sum_{i=1}^{k-1}(-1)^{i+1}\sum_{
|I|\subseteq \{2,\ldots,k\}}\ell(\bigcap_{1<j\in I}(M_j+M_1)/M_1)\\  
&=&\sum_{i=1}^{k-1}(-1)^{i+1} \sum_{|I|\subseteq
\{2,\ldots,k\}}\ell(M_I/(M_1\cap M_I))\\   
&=&\sum_{i=1}^{k-1}(-1)^{i+1}
\sum_{|I|\subseteq \{2,\ldots,k\}}[\ell(M_I)-\ell(M_{I\cup\{1\}})]
\end{eqnarray*}
The terms $\ell(M_I)$ in the last sum 
are all the summands in the sum of the theorem
without the index $1$. The terms $\ell(M_{I\cup\{1\}})$ 
together with $\ell(M_1)$
make up all those who do use the index $1$.
\qed
\end{prop}

\begin{prop}
\mylabel{prop-length}
In the context of Theorem \ref{thm-leykin}, 
let $M_I=R_n\left[\prod_{j\not\in I}{H_j}^{-1}\right]$.
The length of $M=R_n[Q^{-1}]$ is determined recursively as
follows where $H^i_\calA(-)$ is local cohomology with supports in the
ideal $\ideal{H_1,\ldots,H_k}$. 
\begin{itemize}
\item If $H^k_\calA(R_n)=0$ then $\ell(M)=\sum(-1)^i\sum_{|I|=i}\ell(M_I)$.

\item If $H^k_\calA(R_n)\not =0$ then
$\ell(M)=\sum(-1)^i\sum_{|I|=i}\ell(M_I)+1$.
\end{itemize}
This information can be obtained from the intersection lattice.
\proof
In the first case the \v Cech complex shows that $M=\sum_{i=1}^k M_i$
and hence all that needs to be shown is that the two usages of
the symbol $M_I$ here and in Proposition \ref{prop-in-ex} 
agree. In other words, we must show that
\[
R_n\left[\prod_{j\in\{1,\ldots,k\}\setminus I}{H_j}^{-1}\right]\cap 
R_n\left[\prod_{j'\in\{1,\ldots,k\}\setminus I'}{H_{j'}}^{-1}\right] 
=R_n\left[\prod_{j\in\{1,\ldots,k\}\setminus (I\cup
I')}{H_j}^{-1}\right]
\]
for all index sets $I,I'$. 
This, however, is
clear. 

If $H^k_\calA(R_n)\not =0$ then the $H_i$ form a regular sequence and
hence we know that this local cohomology module is of length one,
a suitable generator being annihilated by all $H_i$. 
The formula follows by considering
$\sum_{|I|=1}M_I$ and $0\to \sum_{|I|=1}M_I\to M\to H^k_\calA(R_n)\to 0$.
 \qed
\end{prop}
\begin{rem}
There are substantially more general results by 
{\`A}lvarez-Montaner, Garc{\'{\i}}a-L{\'o}pez,
             and Zarzuela-Armengou. 
In fact, Propositions \ref{prop-in-ex}
and \ref{prop-length} can be modified to apply to the characteristic
cycle of $R_n[Q^{-1}]$. This idea is discussed in \cite{AGZ} and then
used to 
 express the lengths of the modules
$H^r_{\calA}(R_n)$ in terms of Betti numbers obtained from the
intersection lattice (even for subspace arrangements).
\end{rem}

\subsection{Some conjectures}\hfill

We now close with conjectures on the generators of $J(Q^s)$ and
$\ann_{D_n}(Q^{-1})$.
\begin{df}
For a central arrangement $\calA=\{H_1,\ldots,H_k\}$ and 
$Q=H_1\cdots H_k$ we define the ideals $I(\calA)$ and
$I_s(\calA)$ as follows. Let
$H_1,\ldots,H_n$ be linearly independent. 
Choose vector fields  $v_i$ with constant coefficients
such that $v_i\action(H_j)=\delta_{i,j}$ for all $1\le i,j\le
n$.

Form
$P_{i,j}(Q)=\frac{H_iH_j}{H_1\cdots
  H_n}(v_i\action(Q)v_j-v_j\action(Q)v_i)\in D_n$;  
$P_{i,j}(Q)$  kills $Q^s$.  Let 
\[
I(Q)=\left\langle\left\{P_{i,j}(Q')\frac{Q}{Q'}:Q'\divides
  Q\right\}\right\rangle\subseteq D_n. 
\]
We define $I_s(Q)$ recursively. If $\deg(Q)=1$, set $I_s(Q)=I(Q)$.
If
$\deg(Q)>1$, 
\[
I_s(Q)=\left\langle\left\{{Q''}^{s+1}P_{i,j}(Q'){Q''}^{-s}:
Q=Q'Q''\right\}\right\rangle\subseteq D_n[s].
\]
It is apparent that $I_s(Q)$ kills $Q^s$ and $I(Q)$ kills $1/Q$. 
\end{df}
\begin{conj}
For any central arrangement $Q$,
\begin{enumerate}
\item the annihilator $\ann_{D_n}(Q^{-1})$ is $I(Q)+\ideal{\euler+k}$;
\item the annihilator $\ann_{D_n[s]}(Q^s)$ is $I_s(Q)+\ideal{\euler-ks}$. 
\end{enumerate}
\end{conj}
There is certainly a considerable amount of redundancy in these
generators. Particularly for generic arrangements much smaller sets
can be taken. The importance of the conjecture lies perhaps more in the fact
that all operators shown are order one. We make some remarks about
this now.

T.~Torrelli \cite{Torrelli-mero} 
has proved that $\ann(Q^{-1})$ is generated in order one 
for the union of a generic arrangement
with a hyperbolic arrangement. 
A divisor $\Div(f)$ 
on $\CC^n$ is called {\em free} if the module of logarithmic
derivations $\der(\log f)=\{\delta\in\der(R_n):\delta(f)\in\ideal
f\}$ is a locally free $R_n$-module. It
is called 
{\em Koszul-free} if one can choose a basis for the logarithmic derivations
such that their top order parts form a regular sequence in
$\gr_{(0,1)}(D_n)$.  
The complex of logarithmic
differentials $\Omega^\bullet(\log f)$ consists (in the algebraic
case) of those
differential forms $\omega\in\Omega^\bullet(R_n[f^{-1}])$ for which both
$f\omega$ and $f d\omega$ 
are regular forms on $\CC^n$. It is a subcomplex of
$\Omega^\bullet(R_n[f^{-1}])$ and (algebraic) 
Logarithmic Comparison is said to hold if the
inclusion is a quasi-isomorphism. 

Let $\tilde I^{\log f}$ 
be the subideal of $\ann(1/f)$
generated by the order one operators introduced by L.~Narvaez-Macarro
\cite{Ucha-thesis} and put $I^{\log f}=D_n\cdot \der(\log f)$.
F.J~Castro and J.M.~Ucha 
\cite{Castro-Ucha-freedual,Castro-Ucha-comp},
using results and ideas of F.J.~Calderon
\cite{Calderon-ASENS}, proved  
that if $f$ is 
Koszul-free then the map from $\Omega^\bullet(\log f)$
to 
the (holomorphic) de Rham complex of the (holonomic) module $\tilde M^{\log
  f}=D_n/\tilde I^{\log f}$ is a quasi-isomorphism, and
 $\tilde I^{\log f}$ and $I^{\log f}$ are holonomically dual.
Hence if $f$ is Koszul-free and $\tilde M^{\log f}$ regular
holonomic, then 
$\tilde M^{\log f}=R_n[f^{-1}]$ if and only  if (holomorphic) 
Logarithmic Comparison
holds.
In his paper \cite{Torrelli-mero}
Torrelli conjectures that 
if $f$ is reduced (but not necessarily Koszul-free) 
then (holomorphic) Logarithmic Comparison  holds for $f$ 
if and only if $\ann(1/f)=\tilde I^{\log f}$.

Terao conjectured in  \cite{Terao-kyoto} 
that (algebraic) 
Logarithmic Comparison holds for any central arrangement (and
more) and 
there is a proof in the analytic case for free quasi-homogeneous
divisors in \cite{Castro-Narvaez-Mond}.
This can via
Torrelli's conjecture be seen as counterpart to our conjecture. Wiens
and Yuzvinsky 
have proved in \cite{Wiens-Yuzvinsky} 
Terao's conjecture for arrangements in $\CC^{\le 4}$, and all tame
arrangements.

\subsection*{Acknowledgments}\hfill

I would like to thank A.~Leykin for telling me about
Theorem \ref{thm-leykin}, and B.~Sturmfels for his encouragement
regarding this work. My thanks go as well to two unknown referees for
their useful comments, and to MPI and Universit\"at Leipzig for their
hospitality during the final stages of the preparation of this
manuscript.

\bibliography{bib}
\bibliographystyle{abbrv}

\end{document}